\documentclass{amsart}%
\usepackage{amsfonts}
\usepackage{amsmath}
\usepackage{amssymb}
\usepackage{graphicx}%
\newtheorem{theorem}{Theorem}
\theoremstyle{plain}

\newtheorem{definition}{Definition}

\newtheorem{lemma}{Lemma}

\newtheorem{proposition}{Proposition}

\numberwithin{equation}{section}
\begin{document}
\title[Singular polynomials]{Singular polynomials for the symmetric groups}
\author{Charles F. Dunkl}
\address{Department of Mathematics\\
University of Virginia\\
Charlottesville, VA 22904-4137}
\email{cfd5z@virginia.edu}
\urladdr{http://www.people.virginia.edu/\symbol{126}cfd5z/}
\thanks{During the preparation of this paper the author was partially supported by NSF
grant DMS 0100539.}
\date{March 16, 2004}

\begin{abstract}
For certain negative rational numbers $\kappa_{0}$, called singular values,
and associated with the symmetric group $S_{N}$ on $N$ objects, there exist
homogeneous polynomials annihilated by each Dunkl operator when the parameter
$\kappa=\kappa_{0}$. It was shown by the author, de Jeu and Opdam (Trans.
Amer. Math. Soc. 346 (1994), 237-256) \ that the singular values are exactly
the values $-\frac{m}{n}$ with $2\leq n\leq N$, $m=1,2,\ldots$ and $\frac
{m}{n}$ is not an integer. This paper constructs for each pair $\left(
m,n\right)  $ satisfying these conditions an irreducible $S_{N}$-module of
singular polynomials for the singular value $-\frac{m}{n}$. The module is of
isotype $\left(  n-1,\left(  n_{1}-1\right)  ^{l},\rho\right)  $ where
$n_{1}=n/\gcd(m,n)$, ~$\rho=N-\left(  n-1\right)  -l\left(  n_{1}-1\right)  $
and $1\leq\rho\leq n_{1}-1$. The singular polynomials are special cases of
nonsymmetric Jack polynomials. The paper presents some formulae for the action
of Dunkl operators on these polynomials valid in general, and a method for
showing the dependence of poles (in the parameter $\kappa$) on the number of
variables. Murphy elements are used to analyze the representation of $S_{N}$
on irreducible spaces of singular polynomials.

\end{abstract}
\maketitle

\section{Introduction}

We will construct polynomials on $\mathbb{R}^{N}$ which are annihilated by
each Dunkl operator associated with the symmetric group $S_{N},$ acting\ by
permutation of coordinates, when the parameter takes on a singular value
$-\frac{m}{n}$ with $2\leq n\leq N$ and $-\frac{m}{n}\notin\mathbb{Z}$. The
group $S_{N}$ is considered as the finite reflection group of type $A_{N-1}$.
Let $\mathbb{N}_{0}$ denote $\left\{  0,1,2,3,\ldots\right\}  $; for
$\alpha\in\mathbb{N}_{0}^{N}$ (called a ``composition'') let $\left|
\alpha\right|  =\sum_{i=1}^{N}\alpha_{i}$ and define the monomial $x^{\alpha}$
to be $\prod_{i=1}^{N}x_{i}^{\alpha_{i}}$; its degree is $\left|
\alpha\right|  $. The length of a composition is $\ell\left(  \alpha\right)
=\max\left\{  j:\alpha_{j}>0\right\}  $. For $1\leq i\leq N$ let
$\varepsilon\left(  i\right)  \in\mathbb{N}_{0}^{N}$ denote the standard basis
element, that is, $\varepsilon\left(  i\right)  _{j}=\delta_{ij}$. Consider
elements of $S_{N}$ as functions on $\{1,2,\ldots,N\}$ then for $x\in
\mathbb{R}^{N}$ and $w\in S_{N}$ let $\left(  xw\right)  _{i}=x_{w\left(
i\right)  }$ for $1\leq i\leq N$; and extend this action to polynomials by
$wf\left(  x\right)  =f\left(  xw\right)  $. This has the effect that
monomials transform to monomials, $w\left(  x^{\alpha}\right)  =x^{w\alpha}$
where $\left(  w\alpha\right)  _{i}=\alpha_{w^{-1}\left(  i\right)  }$ for
$\alpha\in\mathbb{N}_{0}^{N}$. (Consider $x$ as a row vector, $\alpha$ as a
column vector, and $w$ as a permutation matrix, with $1$'s at the $\left(
w\left(  j\right)  ,j\right)  $ entries.) The reflections in $S_{N}$ are the
transpositions, denoted by $\left(  i,j\right)  $ for $i\neq j$, interchanging
$x_{i}$ and $x_{j}$.

In \cite{D1} the author constructed for each finite reflection group a
parametrized commutative algebra of differential-difference operators. Let
$\kappa$ be a formal parameter, that is, $\mathbb{Q}\left(  \kappa\right)  $
is a transcendental extension of $\mathbb{Q}$. For the symmetric group the
operators are defined as follows:

\begin{definition}
For any polynomial $f$ on $\mathbb{R}^{N}$ and $1\leq i\leq N$ let
\[
\mathcal{D}_{i}f\left(  x\right)  =\frac{\partial}{\partial x_{i}}f\left(
x\right)  +\kappa\sum_{j\neq i}\frac{f\left(  x\right)  -\left(  ij\right)
f\left(  x\right)  }{x_{i}-x_{j}}.
\]

\end{definition}

The polynomials under consideration are elements of $\mathrm{span}%
_{\mathbb{Q}\left(  \kappa\right)  }\left\{  x^{\alpha}:\alpha\in
\mathbb{N}_{0}^{N}\right\}  $. It was shown in \cite{D1} that $\mathcal{D}%
_{i}\mathcal{D}_{j}=\mathcal{D}_{j}\mathcal{D}_{i}$ for $1\leq i,j\leq N$ and
each $\mathcal{D}_{i}$ maps homogeneous polynomials to homogeneous
polynomials. A specific numerical parameter value $\kappa_{0}$ is said to be a
\textit{singular value} (associated with $S_{N}$) if there exists a nonzero
polynomial $p$ such that $\mathcal{D}_{i}p=0$ for $1\leq i\leq N$ when
$\kappa$ is specialized to $\kappa_{0}$, and $p$ is called a \textit{singular
polynomial}. It was shown in \cite{DJO} that the singular values are the
numbers $-\frac{j}{n}$ where $n=2,\ldots,N,\,j\in\mathbb{N}$ and $\frac{j}%
{n}\notin\mathbb{Z}$. The space of homogeneous polynomials of degree $n$,
denoted by $\mathcal{P}_{n}$, is $\mathrm{span}_{\mathbb{Q}\left(
\kappa\right)  }\left\{  x^{\alpha}:\alpha\in\mathbb{N}_{0}^{N},\left\vert
\alpha\right\vert =n\right\}  $. The set of partitions of length $\leq N$ is
denoted by $\mathbb{N}_{0}^{N,P}$ and consists of all $\lambda\in
\mathbb{N}_{0}^{N}$ such that $\lambda_{i}\geq\lambda_{i+1}$ for $1\leq i\leq
N-1$. When writing partitions it is customary to suppress trailing zeros and
to use exponents to indicate multiplicity, for example $\left(  5,2^{3}%
\right)  $ is the same as $\left(  5,2,2,2,0\right)  \in\mathbb{N}_{0}^{5,P}.$
The irreducible representations of $S_{N}$ are labeled by partitions of $N$
(that is, $\tau\in\mathbb{N}_{0}^{N,P}$ and $\left\vert \tau\right\vert =N$)
and we say a polynomial $f$ is of \textit{isotype} $\tau$ if $f$ is an element
of an irreducible $S_{N}$-submodule of $\mathcal{P}_{n}$ on which the
representation $\tau$ is realized. It was conjectured in \cite{DJO} that the
two-part representations $\left(  \mu,N-\mu\right)  $ (with $2\mu\geq N$) give
rise to singular polynomials for the singular values $-\frac{m}{\mu+1}$ with
$\gcd\left(  m,\mu+1\right)  <\frac{\mu+1}{N-\mu}$ (this was shown in
\cite{D2}), and the representations $\left(  s\left(  \mu+1\right)
-1,\mu,\ldots,\mu,\rho\right)  $ for $s,\mu\in\mathbb{N}$ give rise to
singular polynomials for the singular values $-\frac{m}{\mu+1}$ with
$\gcd\left(  m,\mu+1\right)  =1$. The latter is the main topic of this paper.
For example, the singular values $-\frac{m}{6}$ for $N=10$ are associated with
the isotypes $\left(  5,5\right)  $ for $m\equiv1,5\operatorname{mod}6$,
$\left(  5,2,2,1\right)  $ for $m\equiv2,4\operatorname{mod}6$, and $\left(
5,1^{5}\right)  $ for $m\equiv3\operatorname{mod}6$.

In the rest of this introduction we present definitions and key properties of
nonsymmetric Jack polynomials, hook-length products for Ferrers diagrams, and
the fundamental partial order on compositions. Section 2 contains detailed
formulae for the action of $\left\{  \mathcal{D}_{i}\right\}  $ on the
polynomials, with emphasis on the poles. The construction of singular
polynomials is presented in Section 3, and there is a key result on the
absence of certain poles when the number of variables (that is, $N$) is small
enough. Murphy's construction \cite{Mu} of the seminormal representations of
$S_{N}$ is used in Section 4 to analyze the irreducible $S_{N}$-modules
generated by singular polynomials. The conclusion in Section 5 concisely
displays the correspondence between pairs $\left(  m,n\right)  ,2\leq n\leq
N,\frac{m}{n}\notin\mathbb{Z}$ and singular polynomials for $\kappa=-\frac
{m}{n}$, and also considers modules of the specializations of the rational
Cherednik algebra, defined in terms of singular polynomials.

Our construction will be in terms of nonsymmetric Jack polynomials. Since
these have coefficients in $\mathbb{Q}\left(  \kappa\right)  $ with poles at
negative rational values of $\kappa$ it will be important to be precise about
these poles. Any further reference to poles will be with respect to $\kappa$.
The related commutative algebra of self-adjoint operators is generated by
\[
\mathcal{U}_{i}f\left(  x\right)  =\mathcal{D}_{i}x_{i}f\left(  x\right)
-\kappa\sum_{j=1}^{i-1}\left(  j,i\right)  f\left(  x\right)  ,1\leq i\leq N.
\]
(this differs by an additive constant from the notation in \cite[Ch.8]{DX}).
The operators act in a triangular manner on monomials.

\begin{definition}
\label{order}For $\alpha\in\mathbb{N}_{0}^{N}$ let $\alpha^{+}$ denote the
unique partition such that $\alpha^{+}=w\alpha$ for some $w\in S_{N}$. For
$\alpha,\beta\in\mathbb{N}_{0}^{N}$ the partial order $\alpha\succ\beta$
($\alpha$ dominates $\beta$) means that $\alpha\neq\beta$ and $\sum_{i=1}%
^{j}\alpha_{i}\geq\sum_{i=1}^{j}\beta_{i}$ for $1\leq j\leq N$; and
$\alpha\vartriangleright\beta$ means that $\left|  \alpha\right|  =\left|
\beta\right|  $ and either $\alpha^{+}\succ\beta^{+}$ or $\alpha^{+}=\beta
^{+}$ and $\alpha\succ\beta$. The notations $\alpha\succeq\beta$ and
$\alpha\trianglerighteq\beta$ include the case that $\alpha=\beta$.
\end{definition}

Acting on the monomial basis of $\mathcal{P}_{n}$ the operators $\mathcal{U}%
_{i}$ have on-diagonal coefficients involving the following ``rank'' function
on $\mathbb{N}_{0}^{N}$. We denote the cardinality of a set $E$ by $\#E$.

\begin{definition}
\label{rankdef}For $\alpha\in\mathbb{N}_{0}^{N}$ and $1\leq i\leq N$ let
\begin{align*}
r\left(  \alpha,i\right)   &  =\#\left\{  j:\alpha_{j}>\alpha_{i}\right\}
+\#\left\{  j:1\leq j\leq i,\alpha_{j}=\alpha_{i}\right\}  ,\\
\xi_{i}\left(  \alpha\right)   &  =\left(  N-r\left(  \alpha,i\right)
\right)  \kappa+\alpha_{i}+1.
\end{align*}

\end{definition}

Clearly for a fixed $\alpha\in\mathbb{N}_{0}^{N}$ the values $\left\{
r\left(  \alpha,i\right)  :1\leq i\leq N\right\}  $ consist of all of
$\left\{  1,\ldots,N\right\}  $, are independent of trailing zeros (that is,
if $\alpha^{\prime}\in\mathbb{N}_{0}^{M},\alpha_{i}^{\prime}=\alpha_{i}$ for
$1\leq i\leq N$ and $\alpha_{i}^{\prime}=0$ for $N<i\leq M$ then $r\left(
\alpha,i\right)  =r\left(  \alpha^{\prime},i\right)  $ for $1\leq i\leq N$),
and $\alpha\in\mathbb{N}_{0}^{N,P}$ if and only if $r\left(  \alpha,i\right)
=i$ for all $i$. Then (see \cite[p.291]{DX}) $\mathcal{U}_{i}x^{\alpha}%
=\xi_{i}\left(  \alpha\right)  x^{\alpha}+q_{\alpha,i}\left(  x\right)  $
where $q_{\alpha,i}\left(  x\right)  $ is a sum of terms $\pm\kappa x^{\beta}$
with $\alpha\vartriangleright\beta$. The nonsymmetric Jack polynomials are the
simultaneous eigenvectors of $\left\{  \mathcal{U}_{i}:1\leq i\leq N\right\}
,$ well-defined for generic $\kappa$. Opdam \cite[p.83]{O} discovered and
studied them in the wider framework of polynomials associated to
crystallographic root systems. There are two useful normalizations of these
polynomials, one is ``monic in $x$'' the other is ``monic in $p$''. The
$p$-basis is defined by the generating function
\[
\sum_{\alpha\in\mathbb{N}_{0}^{N}}p_{\alpha}\left(  x\right)  y^{\alpha}%
=\prod_{i=1}^{N}\left(  \left(  1-x_{i}y_{i}\right)  ^{-1}\prod_{j=1}%
^{N}\left(  1-x_{i}y_{j}\right)  ^{-\kappa}\right)  ,\text{ for }\max
_{i,j}\left|  x_{i}\right|  \left|  y_{j}\right|  <1.
\]
In contrast to the monomial basis $\mathcal{U}_{i}p_{\alpha}=\xi_{i}\left(
\alpha\right)  p_{\alpha}+q_{\alpha,i}^{\prime}$ where $q_{\alpha,i}^{\prime}$
is a sum of terms $\pm\kappa p_{\beta}$ with $\beta\vartriangleright\alpha$
(and $\ell\left(  \beta\right)  \leq\ell\left(  \alpha\right)  $), (see
\cite[Prop. 8.4.11]{DX}).

\begin{definition}
For $\alpha\in\mathbb{N}_{0}^{N}$ let $\zeta_{\alpha},\zeta_{\alpha}^{x}$
denote the $p$-monic and $x$-monic, respectively, simultaneous eigenvectors,
that is, $\mathcal{U}_{i}\zeta_{\alpha}=\xi_{i}\left(  \alpha\right)
\zeta_{\alpha},\,\mathcal{U}_{i}\zeta_{\alpha}^{x}=\xi_{i}\left(
\alpha\right)  \zeta_{\alpha}^{x}$ for $1\leq i\leq N$ and $\zeta_{\alpha
}=p_{\alpha}+\sum_{\beta\vartriangleright\alpha}A_{\beta\alpha}p_{\beta
},\,\zeta_{\alpha}^{x}=x^{\alpha}+\sum_{\alpha\vartriangleright\beta}%
A_{\beta\alpha}^{x}x^{\beta}$, with coefficients $A_{\beta\alpha}%
,A_{\beta\alpha}^{x}\in\mathbb{Q}\left(  \kappa\right)  $.
\end{definition}

Suppose that $\ell\left(  \alpha\right)  =m$ for some $m\geq1$ then the
coefficients $A_{\beta\alpha}$ do not depend on $N\geq m$ (and $A_{\beta
\alpha}\neq0$ implies $\ell\left(  \beta\right)  \leq m$); on the other hand,
if $\beta\in\mathbb{N}_{0}^{N}$ and $\ell\left(  \beta\right)  \leq m $ then
$A_{\beta\alpha}^{x}$ does not depend on $N\geq m$ (that is, if $N>M\geq m$
then the projection of $\mathbb{R}^{N}$ onto $\mathbb{R}^{M}$ setting
$x_{M+1}=\ldots=x_{N}=0$ and annihilating the terms $A_{\beta\alpha}%
^{x}x^{\beta}$ with $\beta_{i}\neq0$ for some $i>M$ ($\ell\left(
\beta\right)  >m$), produces the $\mathbb{R}^{M}$-polynomial. The relation
between the two types involves hook-length products. Suppose $\lambda
\in\mathbb{N}_{0}^{N,P}$ and $\ell\left(  \lambda\right)  =m$; the Ferrers
diagram of $\lambda$ is the set $\left\{  \left(  i,j\right)  :1\leq i\leq
m,1\leq j\leq\lambda_{i}\right\}  .$ Each node $\left(  i,j\right)  $ has the
\textit{arm} $\left\{  \left(  i,l\right)  :j<l\leq\lambda_{i}\right\}  $ and
the \textit{leg} $\left\{  \left(  l,j\right)  :i<l,j\leq\lambda_{l}\right\}
$. The node itself, the arm and the leg make up the \textit{hook}. For
$t\in\mathbb{Q}\left(  \kappa\right)  $ the \textit{hook-length}, the
hook-length product and generalized Pochhammer symbol for $\lambda$ are given
by
\begin{align*}
h\left(  \lambda,t;i,j\right)   &  =\left(  \lambda_{i}-j+t+\kappa\#\left\{
l:i<j,\,j\leq\lambda_{l}\right\}  \right) \\
h\left(  \lambda,t\right)   &  =\prod_{i=1}^{m}\prod_{j=1}^{\lambda_{i}%
}h\left(  \lambda,t;i,j\right)  ,\\
\left(  t\right)  _{\lambda}  &  =\prod_{i=1}^{m}\prod_{j=1}^{\lambda_{i}%
}\left(  t-(i-1)\kappa+j-1\right)  .
\end{align*}
The coordinate-wise notation for hook-lengths will appear in the context of
specializations of $\kappa$ to negative rational numbers. The compositions
$\alpha\in\mathbb{N}_{0}^{N}$ are associated with the products
\[
\mathcal{E}_{\varepsilon}\left(  \alpha\right)  =\prod\left\{  1+\frac
{\varepsilon\kappa}{\kappa\left(  r\left(  \alpha,i\right)  -r\left(
\alpha,j\right)  \right)  +\alpha_{j}-\alpha_{i}}:i<j,\,\alpha_{i}<\alpha
_{j}\right\}  ,\varepsilon=\pm.
\]
Note that the denominator is identical to $\xi_{j}\left(  \alpha\right)
-\xi_{i}\left(  \alpha\right)  $ and $\mathcal{E}_{\varepsilon}\left(
\lambda\right)  =1$ for $\lambda\in\mathbb{N}_{0}^{N,P}$. Then (see
\cite[p.323]{DX}) for each $\alpha\in\mathbb{N}_{0}^{N}$%
\[
\zeta_{\alpha}=\mathcal{E}_{+}\left(  \alpha\right)  \mathcal{E}_{-}\left(
\alpha\right)  \frac{h\left(  \alpha^{+},\kappa+1\right)  }{h\left(
\alpha^{+},1\right)  }\zeta_{\alpha}^{x}.
\]
Also $\zeta_{\alpha}\left(  1^{N}\right)  =\mathcal{E}_{-}\left(
\alpha\right)  \left(  N\kappa+1\right)  _{\alpha^{+}}\,/h\left(  \alpha
^{+},1\right)  $. Knop and Sahi \cite{KS} by finding explicit combinatorial
formulae in terms of tableaux established the key theorem that $h\left(
\lambda,\kappa+1\right)  \zeta_{\lambda}^{x}$ is a polynomial with
coefficients in $\mathbb{Z}\left[  \kappa\right]  $. However as $N$ decreases
the set of $\kappa$-poles of $\zeta_{\lambda}^{x}$ also decreases, and
specific results will be established and used in the sequel (by \cite[Cor.
4.7]{KS} the coefficient of $x_{m+1}\ldots x_{m+n}$ in $\zeta_{\lambda}^{x} $
is $n!\,\kappa^{n}/h\left(  \lambda,\kappa+1\right)  $ where $\left\vert
\lambda\right\vert =n$ and $\ell\left(  \lambda\right)  =m$; thus for $m\leq
N<m+n$ one expects some poles to be omitted). An obvious sufficient condition
for the presence of a pole in $\zeta_{\lambda}^{x}$ for given $N$ is its
presence in $\zeta_{\lambda}^{x}\left(  1^{N}\right)  =\left(  N\kappa
+1\right)  _{\lambda}\,/h\left(  \lambda,\kappa+1\right)  .$

We can now state our main results: for each isotype $\tau$ and singular value
$\kappa_{0}$ the corresponding singular polynomials form the $S_{N}$-module
generated by $\zeta_{\lambda}^{x}$ for a certain $\lambda$, that is,
$\mathrm{span}_{\mathbb{Q}}\left\{  w\zeta_{\lambda}^{x}:w\in S_{N}\right\}
$; (in fact a basis will be specified in terms of reverse lattice permutations
of $\lambda$)

\begin{itemize}
\item for $\tau=\left(  \mu,N-\mu\right)  ,\kappa_{0}=-\frac{m}{\mu+1}$ with
$\gcd\left(  m,\mu+1\right)  <\frac{\mu+1}{N-\mu}$, let $\lambda=\left(
m^{N-\mu},0^{\mu}\right)  $ (that is, $m$ is repeated $N-\mu$ times followed
by $\mu$ zeros)

\item for $\tau=\left(  s\left(  \mu+1\right)  +\mu,\mu^{l},\rho\right)  $
where $l\geq1,s\geq0,$ and $1\leq\rho\leq\mu$ (so that $N=\left(
s+l+1\right)  \mu+s+\rho$), $\kappa_{0}=-\frac{m}{\mu+1}$ with $\gcd\left(
m,\mu+1\right)  =1$, let $\lambda=\left(  \left(  m\left(  s+l+1\right)
\right)  ^{\rho},\left(  m\left(  s+l\right)  \right)  ^{\mu},\ldots,\left(
m\left(  s+1\right)  \right)  ^{\mu},0^{s\left(  \mu+1\right)  +\mu}\right)  $.
\end{itemize}

For example, let $N=10,\tau=\left(  5,2,2,1\right)  ,\kappa_{0}=-\frac{m}{3}$
and $\gcd\left(  m,3\right)  =1$, then $\lambda=\left(  4m,3m,3m,2m,2m,0^{5}%
\right)  $. The singular polynomials for $N=2k+1,\,\tau=\left(
2k-1,1,1\right)  ,\,\kappa_{0}=-\frac{m}{2},\,\lambda=\left(  m\left(
k+1\right)  ,mk\right)  $ were found in \cite{D3} by a different method (not
suitable for the general problem).

\section{Differentiation Formulae}

This section contains expressions for $\mathcal{D}_{i}\zeta_{\alpha}$ in terms
of $\left\{  \zeta_{\beta}:\left|  \beta\right|  =\left|  \alpha\right|
-1\right\}  $, valid for generic $\kappa$. There is some material dealing with
$x$-monic polynomials, however the $p$-monic polynomials have somewhat simpler
formulae. The basic step is the formula for $\mathcal{D}_{m}\zeta_{\alpha}$
where $\ell\left(  \alpha\right)  =m,$ $\alpha\in\mathbb{N}_{0}^{N}$ $;$
further from properties of the $p$-basis it follows that $\mathcal{D}_{i}%
\zeta_{\alpha}=0$ for $i>m$. The computation involves a cyclic shift. For
$1\leq i\leq j\leq N$ let $\left[  i,j\right]  $ denote the interval $\left\{
k\in\mathbb{Z}:i\leq k\leq j\right\}  $ and let $S\left[  i,j\right]  $ denote
the subgroup of $S_{N}$ generated by $\left\{  \left(  i,i+1\right)
,(i+1,i+2),\ldots,\left(  j-1,j\right)  \right\}  $ (isomorphic to $S_{j+1-i}$).

\begin{definition}
For $1<m\leq N$ let $\theta_{m}=\left(  1,2\right)  \left(  2,3\right)
\ldots\left(  m-1,m\right)  \in S_{N}$, and if $\alpha\in\mathbb{N}_{0}^{N}$
satisfies $\ell\left(  \alpha\right)  =m$ then $\widetilde{\alpha}=\theta
_{m}\left(  \alpha-\varepsilon\left(  m\right)  \right)  =\left(  \alpha
_{m}-1,\alpha_{1},\ldots,\alpha_{m-1},0,\ldots\right)  $.
\end{definition}

\begin{lemma}
\label{cyclic}Suppose $\alpha\in\mathbb{N}_{0}^{N}$ satisfies $\ell\left(
\alpha\right)  =m$ then (i) $\mathcal{U}_{1}\theta_{m}\mathcal{D}_{m}%
\zeta_{\alpha}=$\newline$\left(  \xi_{m}\left(  \alpha\right)  -1\right)
\theta_{m}\mathcal{D}_{m}\zeta_{\alpha}$, (ii) $\mathcal{U}_{i}\theta
_{m}\mathcal{D}_{m}\zeta_{\alpha}=\xi_{i-1}\left(  \alpha\right)  \theta
_{m}\mathcal{D}_{m}\zeta_{\alpha}$ for $1<i\leq m$, and (iii) $\mathcal{U}%
_{i}\theta_{m}\mathcal{D}_{m}\zeta_{\alpha}=\left(  \left(  N-i\right)
\kappa+1\right)  \theta_{m}\mathcal{D}_{m}\zeta_{\alpha}$ for $i>m$.
\end{lemma}

\begin{proof}
The commutation $\left(  x_{m}\mathcal{D}_{m}-\mathcal{D}_{m}x_{m}\right)
f=-f-\kappa\sum_{j\neq m}\left(  j,m\right)  f$ (\cite[p.290]{DX}) shows that
\begin{align*}
\mathcal{D}_{m}x_{m}\mathcal{D}_{m}\zeta_{\alpha}  &  =-\mathcal{D}_{m}%
\zeta_{\alpha}+\mathcal{D}_{m}\left(  \mathcal{D}_{m}x_{m}-\kappa\sum
_{j<m}\left(  j,m\right)  \right)  \zeta_{\alpha}-\kappa\sum_{j>m}%
\mathcal{D}_{m}\left(  j,m\right)  \zeta_{\alpha}\\
&  =\mathcal{D}_{m}\left(  \xi_{m}\left(  \alpha\right)  -1\right)
\zeta_{\alpha}%
\end{align*}
because $\mathcal{D}_{m}\left(  j,m\right)  \zeta_{\alpha}=\left(  j,m\right)
\mathcal{D}_{j}\zeta_{\alpha}=0$ for $j>m$. Apply $\theta_{m}$ to the previous
equation to prove part (i) (since $\theta_{m}\mathcal{D}_{m}x_{m}%
=\mathcal{D}_{1}x_{1}\theta_{m}$). Next suppose that $1<i\leq m$ then
$\theta_{m}^{-1}\mathcal{U}_{i}\theta_{m}=\mathcal{D}_{i-1}x_{i-1}-\kappa
\sum_{j=1}^{i-2}\left(  j,i-1\right)  -\kappa\left(  m,i-1\right)  $. Apply
this operator to $\mathcal{D}_{m}\zeta_{\alpha}$ to obtain
\begin{align*}
\theta_{m}^{-1}\mathcal{U}_{i}\theta_{m}\mathcal{D}_{m}\zeta_{\alpha}  &
=\mathcal{D}_{m}\left(  \mathcal{D}_{i-1}x_{i-1}-\kappa\sum_{j=1}^{i-2}\left(
j,i-1\right)  \right)  \zeta_{\alpha}\\
&  +\kappa\left(  \mathcal{D}_{i-1}\left(  i-1,m\right)  -\left(
m,i-1\right)  \mathcal{D}_{m}\right) \\
&  =\mathcal{D}_{m}\xi_{i-1}\left(  \alpha\right)  \zeta_{\alpha}.
\end{align*}
The computation uses the commutativity of $\mathcal{D}_{m}$ and $\mathcal{D}%
_{i-1}$ and the commutation $\left(  x_{j}\mathcal{D}_{m}-\mathcal{D}_{m}%
x_{j}\right)  f=\kappa\left(  j,m\right)  f$ (\cite[p.290]{DX}) for $j\neq
m$). This shows part (ii). Similarly for $i>m$ we have that $\theta_{m}%
^{-1}\mathcal{U}_{i}\theta_{m}=\mathcal{U}_{i}$ and $\mathcal{U}%
_{i}\mathcal{D}_{m}-\mathcal{D}_{m}\mathcal{U}_{i}=\kappa\left(
\mathcal{D}_{i}\left(  i,m\right)  -\left(  i,m\right)  \mathcal{D}%
_{m}+\mathcal{D}_{m}\left(  i,m\right)  \right)  =\kappa\left(  i,m\right)
\mathcal{D}_{i}$. But $\mathcal{D}_{i}\zeta_{\alpha}=0$ for $i>m$ and so
$\mathcal{U}_{i}\mathcal{D}_{m}\zeta_{\alpha}=\xi_{i}\left(  \alpha\right)
\mathcal{D}_{m}\zeta_{\alpha}$ and $\xi_{i}\left(  \alpha\right)  =\left(
N-i\right)  \kappa+1$; proving part (iii).
\end{proof}

The following is used to pick out a coefficient in $\mathcal{D}_{m}%
\zeta_{\alpha}$.

\begin{lemma}
\label{difp}Suppose $\alpha,\beta\in\mathbb{N}_{0}^{N},$ $\left|
\alpha\right|  =\left|  \beta\right|  $ and $\ell\left(  \alpha\right)
=\ell\left(  \beta\right)  =m$, if $p_{\alpha-\varepsilon\left(  m\right)  }$
appears with a nonzero coefficient in the expansion of $\mathcal{D}%
_{m}p_{\beta}$ then either $\alpha=\beta$ or $\beta\vartriangleleft\alpha$. If
$\alpha=\beta$ then the coefficient is $\left(  N+1-r\left(  \alpha,m\right)
\right)  \kappa+\alpha_{m}$.
\end{lemma}

\begin{proof}
By Prop.8.4.3 \cite[p.294]{DX}
\begin{gather*}
\mathcal{D}_{m}p_{\beta}=\left(  \left(  N-\#\left\{  j:\beta_{j}\geq\beta
_{m}\right\}  +1\right)  \kappa+\beta_{m}\right)  p_{\beta-\varepsilon\left(
m\right)  }\\
+\kappa\sum\left\{  p_{\gamma}:\gamma=\beta+n\varepsilon\left(  m\right)
-\left(  n+1\right)  \varepsilon\left(  j\right)  ,\max\left(  0,\beta
_{j}-\beta_{m}\right)  \leq n\leq\beta_{j}-1,j\neq m\right\} \\
-\kappa\sum\left\{  p_{\gamma}:\gamma=\beta-\left(  n+1\right)  \varepsilon
\left(  m\right)  +n\varepsilon\left(  j\right)  ,\max\left(  1,\beta
_{m}-\beta_{j}\right)  \leq n\leq\beta_{m}-1,j\neq m\right\}  .
\end{gather*}
If $\beta=\alpha$ then the coefficient of $p_{\alpha-\varepsilon_{m}}$ is
$\left(  N-r\left(  \alpha,m\right)  +1\right)  \kappa+\alpha_{m}$; note that
$j>m$ implies $\alpha_{j}=0<\alpha_{m}$. If $p_{\alpha-\varepsilon\left(
m\right)  }$ has the coefficient $\kappa$ then $\beta=\alpha-\left(
n+1\right)  \left(  \varepsilon\left(  m\right)  -\varepsilon\left(  j\right)
\right)  $ and the restrictions are equivalent to $0\leq n\leq\alpha
_{m}-\alpha_{j}-2$ (for some $j\neq m$) ; thus $n+1<\alpha_{m}-\alpha_{j}$ and
$\beta^{+}\prec\alpha^{+}$ by Lemma 8.2.3(iv) \cite[p.289]{DX}. If
$p_{\alpha-\varepsilon\left(  m\right)  }$ has the coefficient $-\kappa$ then
$\beta=\alpha-n\left(  \varepsilon\left(  j\right)  -\varepsilon\left(
m\right)  \right)  $ and the restrictions are equivalent to $1\leq n\leq
\alpha_{j}-\alpha_{m}$ (for some $j\neq m$). If $n<\alpha_{j}-\alpha_{m}$ then
by the same lemma $\beta^{+}\prec\alpha^{+}$. If $n=\alpha_{j}-\alpha_{m}$
then $\beta=\left(  j,m\right)  \alpha$ and $\beta\prec\alpha$ because
$\alpha_{j}>\alpha_{m}$ and $j<m$ (using the hypothesis $\ell\left(
\alpha\right)  =m$ ). Thus, if $\alpha\neq\beta$ then $\beta\vartriangleleft
\alpha$.
\end{proof}

\begin{theorem}
Suppose $\alpha\in\mathbb{N}_{0}^{N}$ and $\ell\left(  \alpha\right)  =m$,
then%
\[
\mathcal{D}_{m}\zeta_{\alpha}=\left(  \left(  N+1-r\left(  \alpha,m\right)
\right)  \kappa+\alpha_{m}\right)  \theta_{m}^{-1}\zeta_{\widetilde{\alpha}}.
\]

\end{theorem}

\begin{proof}
By Lemma \ref{cyclic}, $\theta_{m}\mathcal{D}_{m}\zeta_{\alpha}$ is a
simultaneous eigenvector of $\left\{  \mathcal{U}_{i}:1\leq i\leq N\right\}  $
with eigenvalues $\left(  \xi_{m}\left(  \alpha\right)  -1,\xi_{1}\left(
\alpha\right)  ,\ldots,\xi_{m-1}\left(  \alpha\right)  ,\xi_{m+1}\left(
\alpha\right)  ,\ldots\right)  $. We claim these are the eigenvalues for
$\widetilde{\alpha}$. Indeed $r\left(  \widetilde{\alpha},1\right)
=\#\left\{  j:\alpha_{j}>\alpha_{m}-1,j<m\right\}  +1=\#\left\{  j:\alpha
_{j}\geq\alpha_{m}\right\}  =r\left(  \alpha,m\right)  .$ For any values
$\alpha_{i}$ different from $\alpha_{m}-1$ and 0 it is obvious that $r\left(
\alpha,i\right)  =r\left(  \widetilde{\alpha},i+1\right)  $. Suppose for some
$i<m$ that $\alpha_{i}=\alpha_{m}-1$ then
\begin{align*}
r\left(  \widetilde{\alpha},i+1\right)   &  =\#\left\{  j:\alpha_{j}%
>\alpha_{m}-1,j<m\right\}  +\#\left\{  j:j\leq i,\alpha_{j}=\alpha
_{m}-1\right\}  +1\\
&  =\#\left\{  j:\alpha_{j}>\alpha_{m}-1,j\leq m\right\}  +\#\left\{  j:j\leq
i,\alpha_{j}=\alpha_{m}-1\right\} \\
&  =r\left(  \alpha,i\right)  .
\end{align*}
For $i>m,$obviously $\xi_{i}\left(  \alpha\right)  =\xi_{i}\left(
\widetilde{\alpha}\right)  =\left(  N-i\right)  \kappa+1$. Thus $\theta
_{m}\mathcal{D}_{m}\zeta_{\alpha}=c\zeta_{\widetilde{\alpha}}$ for some
constant $c$, which will be determined by finding the coefficient of
$\theta_{m}^{-1}p_{\widetilde{\alpha}}=p_{\alpha-\varepsilon\left(  m\right)
}$ in $\mathcal{D}_{m}\zeta_{\alpha}$. Since $\zeta_{\alpha}=p_{\alpha}%
+\sum_{\beta\vartriangleright\alpha}A_{\beta\alpha}p_{\beta}$ (and
$\ell\left(  \beta\right)  \leq m$) we obtain $\mathcal{D}_{m}\zeta_{\alpha
}=\mathcal{D}_{m}p_{\alpha}+\sum_{\beta\vartriangleright\alpha}A_{\beta\alpha
}\mathcal{D}_{m}p_{\beta}$. If $p_{\alpha-\varepsilon\left(  m\right)  }$ has
a nonzero coefficient in $\mathcal{D}_{m}p_{\beta}$ then $\ell\left(
\beta\right)  =m$ (else $\beta_{m}=0$ and $\mathcal{D}_{m}p_{\beta}=0$) and by
Lemma \ref{difp} $\beta=\alpha$ or $\beta\vartriangleright\alpha$. only the
case $\beta=\alpha$ can occur and thus $c=\left(  N+1-r\left(  \alpha
,m\right)  \right)  \kappa+\alpha_{m}$.
\end{proof}

With the intention of using the theorem to compute arbitrary $\mathcal{D}%
_{i}\zeta_{\lambda}$ with $\lambda\in\mathbb{N}_{0}^{N,P}$ we observe that it
suffices to consider the points of decrease, that is, $\lambda_{i}%
>\lambda_{i+1},$(the values of $i$ for which $\lambda-\varepsilon\left(
i\right)  $ is a partition) then apply the transpositions $\left(
j,j+1\right)  $ successively for $j=i,i+1,\ldots,\ell\left(  \lambda\right)
-1,$ apply $\mathcal{D}_{m}$, with a result involving $\zeta_{\alpha}$ where
$\alpha=\left(  \lambda_{i}-1,\lambda_{1},\ldots,\lambda_{i-1},\lambda
_{i+1,\ldots}\right)  $ (this is an over-simplification; actually all the
points of decrease between $i$ and $\ell\left(  \lambda\right)  $ must be
considered). Finally transform $\zeta_{\lambda-\varepsilon\left(  i\right)  }$
to $\zeta_{\alpha}$ with another sequence of transpositions. As mentioned
before it is necessary to keep track of the $\kappa$-poles occurring in these
operations. The basic step is the action of an adjacent transposition on
$\zeta_{\alpha}.$

\begin{proposition}
\label{z2sz}Suppose $\alpha\in\mathbb{N}_{0}^{N}$, and $\alpha_{i}%
>\alpha_{i+1}$ for some $i$, then let $\sigma=\left(  i,i+1\right)  $ and
$a=\kappa\left(  \left(  r\left(  \alpha,i+1\right)  -r\left(  \alpha
,i\right)  \right)  \kappa+\alpha_{i}-\alpha_{i+1}\right)  ^{-1}$ then
$\zeta_{\sigma\alpha}=\sigma\zeta_{\alpha}-a\zeta_{\alpha}$ and $\zeta
_{\sigma\alpha}^{x}=\dfrac{1}{1-a^{2}}\left(  \sigma\zeta_{\alpha}^{x}%
-a\zeta_{\alpha}^{x}\right)  $.
\end{proposition}

The proof for the $p$-monic case is in Prop. 8.5.5 \cite[p.301]{DX}; the proof
for the $x$-monic case can be deduced from the inverse of the $p$-monic
formula and the equation $\zeta_{\alpha}^{x}=\sigma\zeta_{\sigma\alpha}%
^{x}+a\zeta_{\sigma\alpha}^{x}$ (arguing that $x^{\alpha}$ does not appear in
$\zeta_{\sigma\alpha}^{x}$ since $\alpha\vartriangleright\sigma\alpha$). Note
that the denominator $\left(  r\left(  \alpha,i+1\right)  -r\left(
\alpha,i\right)  \right)  \kappa+\alpha_{i}-\alpha_{i+1}=\xi_{i}\left(
\alpha\right)  -\xi_{i+1}\left(  \alpha\right)  $. For singular values of
$\kappa$ it can happen that $a=-1$ implying that $\sigma\zeta_{\alpha}%
^{x}=-\zeta_{\alpha}^{x}$. We need an extension of the proposition applying to
the situation of several adjacent entries of $\alpha$ being equal.

\begin{proposition}
\label{movert}Suppose $\alpha\in\mathbb{N}_{0}^{N}$ with $\alpha
_{i}=a>b=\alpha_{i+j}$ for $1\leq j\leq s$ , where $1\leq i<i+s\leq N$, then
\[
\zeta_{\left(  i,i+s\right)  \alpha}=\left(  i,i+s\right)  \zeta_{\alpha
}-\frac{\kappa}{\left(  r\left(  \alpha,i+s\right)  -r\left(  \alpha,i\right)
\right)  \kappa+a-b}\left(  1+\sum_{j=1}^{s-1}\left(  i,i+j\right)  \right)
\zeta_{\alpha}.
\]

\end{proposition}

\begin{proof}
Observe that $\left(  i,i+s\right)  \alpha=\left(  \ldots,b,\ldots
,b,a,\ldots\right)  $. The proof is by induction on $s$ and depends on the
invariance of $\zeta_{\alpha}$ under the subgroup $S\left[  i+1,i+s\right]  $.
Since $r\left(  \alpha,i+j\right)  =r\left(  \alpha,i+1\right)  +j-1$ for
$1\leq j\leq s$ let $C=\left(  r\left(  \alpha,i+1\right)  -r\left(
\alpha,i\right)  -1\right)  \kappa+a-b$ so that $\left(  r\left(
\alpha,i+j\right)  -r\left(  \alpha,i\right)  \right)  \kappa+a-b=j\kappa+C$ ,
also let $f_{j}=\zeta_{\left(  i,i+j\right)  \alpha}$ and $c_{j}=-\frac
{\kappa}{j\kappa+C}$ . By Proposition \ref{z2sz} $f_{j+1}=\left(
i+j,i+j+1\right)  f_{j}+c_{j+1}f_{j}.$ By the inductive hypothesis
$f_{j}=\left(  i,i+j\right)  \zeta_{\alpha}+c_{j}\left(  1+\sum_{k=1}%
^{j-1}\left(  i,i+k\right)  \right)  \zeta_{\alpha}.$Thus
\begin{align*}
f_{j+1}  &  =\left(  i+j,i+j+1\right)  \left(  i,i+j\right)  \zeta_{\alpha
}+c_{j+1}\left(  i,i+j\right)  \zeta_{\alpha}\\
&  +c_{j}\left(  \left(  i+j,i+j+1\right)  +c_{j+1}\right)  \left(
1+\sum_{k=1}^{j-1}\left(  i,i+k\right)  \right)  \zeta_{\alpha}\\
&  =\left(  i,i+j\right)  \zeta_{\alpha}+c_{j+1}\left(  i,i+j\right)
\zeta_{\alpha}\\
&  +c_{j}\left(  1+\sum_{k=1}^{j-1}\left(  i,i+k\right)  \right)  \left(
\left(  i+j,i+j+1\right)  +c_{j+1}\right)  \zeta_{\alpha}\\
&  =\left(  i,i+j\right)  \zeta_{\alpha}+c_{j+1}\left(  i,i+j\right)
\zeta_{\alpha}+c_{j}\left(  1+c_{j+1}\right)  \left(  1+\sum_{k=1}%
^{j-1}\left(  i,i+k\right)  \right)  \zeta_{\alpha}.
\end{align*}
By the invariance property of $\zeta_{\alpha}$ we have $\left(
i+j,i+j+1\right)  \left(  i,i+j\right)  \zeta_{\alpha}=$\newline$\left(
i+j,i+j+1\right)  \left(  i,i+j\right)  \left(  i+j,i+j+1\right)
\zeta_{\alpha}=\allowbreak\left(  i,i+j\right)  \zeta_{\alpha}$. Furthermore
\newline$c_{j}\left(  1+c_{j+1}\right)  =-\frac{\kappa}{j\kappa+C}\left(
1-\frac{\kappa}{j\kappa+\kappa+C}\right)  =c_{j+1}$ and this completes the induction.
\end{proof}

There is a similar result for the opposite direction.

\begin{proposition}
\label{movelt}Suppose $\alpha\in\mathbb{N}_{0}^{N}$ with $\alpha
_{i+j}=b>a=\alpha_{i+s}$ for $0\leq j\leq s-1$ , where $1\leq i<i+s\leq N$,
then
\[
\zeta_{\left(  i,i+s\right)  \alpha}=\left(  i,i+s\right)  \zeta_{\alpha
}-\frac{\kappa}{\left(  r\left(  \alpha,i+s\right)  -r\left(  \alpha,i\right)
\right)  \kappa+b-a}\left(  1+\sum_{j=1}^{s-1}\left(  i+j,i+s\right)  \right)
\zeta_{\alpha}.
\]

\end{proposition}

\begin{proof}
Proceeding similarly to the previous case, for $1\leq j\leq s$, $r\left(
\alpha,i+s-j\right)  =r\left(  \alpha,i+s-1\right)  +1-j$ and let $f_{j}%
=\zeta_{\left(  i+s-j,i+s\right)  \alpha}$ and \newline$c_{j}=-\kappa\left(
\left(  r\left(  \alpha,i+s\right)  -r\left(  \alpha,i+s-1\right)
-1+j\right)  \kappa+b-a\right)  ^{-1}$. Also $\zeta_{\alpha}$ is invariant
under $S\left[  i,i+s-1\right]  $. The inductive step is based on\newline%
$f_{j+1}=\left(  i+s-j-1,i+s-j\right)  f_{j}+c_{j+1}f_{j}$ (and $f_{0}%
=\zeta_{\alpha}$). The rest of the argument is similar to the previous one and
is omitted.
\end{proof}

To illustrate the basic step, apply $\mathcal{D}_{i+s}$ to both sides of the
formula in Proposition \ref{movert} and obtain
\begin{gather*}
\mathcal{D}_{i+s}\zeta_{\left(  i,i+s\right)  \alpha}=\left(  i,i+s\right)
\mathcal{D}_{i}\zeta_{\alpha}\\
-\frac{\kappa}{\left(  r\left(  \alpha,i+s\right)  -r\left(  \alpha,i\right)
\right)  \kappa+a-b}\left(  1+\sum_{j=1}^{s-1}\left(  i,i+j\right)  \right)
\mathcal{D}_{i+s}\zeta_{\alpha},\\
\mathcal{D}_{i}\zeta_{\alpha}=\left(  i,i+s\right)  \mathcal{D}_{i+s}%
\zeta_{\left(  i,i+s\right)  \alpha}\\
+\frac{\kappa}{\left(  r\left(  \alpha,i+s\right)  -r\left(  \alpha,i\right)
\right)  \kappa+a-b}\left(  i,i+s\right)  \left(  1+\sum_{j=1}^{s-1}\left(
i,i+j\right)  \right)  \mathcal{D}_{i+s}\zeta_{\alpha};
\end{gather*}
so that the index for $\mathcal{D}$ is increased (eventually to $m$).

Fix a partition $\lambda\in\mathbb{N}_{0}^{N,P}$, suppose that the parts of
$\lambda$ have $M$ distinct nonzero values, the points of decrease are
$i_{1}<i_{2}<\ldots<i_{M}$, so that $\lambda_{i}$ is constant on each interval
$i_{j-1}<i\leq i_{j}$ \ (interpret $i_{0}=0$, also let $\ell\left(
\lambda\right)  =m=i_{M}$). For $1\leq j<k\leq M$ let
\begin{align*}
C_{jk}  &  =\frac{\kappa}{\left(  i_{k}-i_{j}\right)  \kappa+\lambda_{i_{j}%
}-\lambda_{i_{k}}},\\
w_{j}  &  =1+\sum_{r=i_{j}+1}^{i_{j+1}-1}\left(  i_{j},r\right)  \in
\mathbb{Z}S\left[  i_{j},i_{j+1}-1\right]  ,\\
z_{jk}  &  =\left(  i_{k-1},i_{k}\right)  -C_{jk}w_{k-1},
\end{align*}
further let $\mu\left(  j,k\right)  \in\mathbb{N}_{0}^{N}$ be the action on
$\lambda$ by the cyclic shift on the interval $\left\{  i_{j},\ldots
,i_{k}\right\}  $, that is $\mu\left(  j,k\right)  _{i_{k}}=\lambda_{i_{j}%
},\,\mu\left(  j,k\right)  _{i}=\lambda_{i+1}$ for $i_{j}\leq i<i_{k}$ and
$\mu\left(  j,k\right)  _{i}=\lambda_{i}$ for $i<i_{j}$ or $i>i_{k}$.
Proposition \ref{movert} applies to the transformation of $\zeta_{\mu\left(
j,k\right)  }$ to $\zeta_{\mu\left(  j,k+1\right)  }$; note that $r\left(
\mu\left(  j,k\right)  ,i_{k}\right)  =i_{j}$ and $r\left(  \mu\left(
j,k\right)  ,i_{k+1}\right)  =i_{k+1}$ thus $\zeta_{\mu\left(  j,k+1\right)
}=z_{j,k+1}\,\zeta_{\mu\left(  j,k\right)  }$. The start of this recurrence is
$\zeta_{\mu\left(  j,j\right)  }=\zeta_{\lambda}$. The object is to express
any $\mathcal{D}_{i}\zeta_{\lambda}$ in terms of $\mathcal{D}_{m}\zeta
_{\mu\left(  j,M\right)  }$, $j=1,\ldots,M$. It suffices to consider $\left\{
\mathcal{D}_{i_{j}}\zeta_{\lambda}\right\}  $ since $\mathcal{D}_{i}%
\zeta_{\lambda}=\left(  i,i_{j}\right)  \mathcal{D}_{i_{j}}\zeta_{\lambda}$
for $i_{j-1}<i<i_{j}$.

\begin{lemma}
For $k=1,\ldots,M-j$
\begin{gather*}
\mathcal{D}_{i_{M}}z_{j,M}\ldots z_{j,M-k+1}=\left(  i_{M-1},i_{M}\right)
\left(  i_{M-2},i_{M-1}\right)  \ldots\left(  i_{M-k},i_{M-k+1}\right)
\mathcal{D}_{i_{M-k}}\\
-\sum_{s=0}^{k-1}C_{j,M-s}\left(  i_{M-1},i_{M}\right)  \ldots\left(
i_{M-s},i_{M-s+1}\right)  w_{M-s-1}z_{j,M-s-1}\ldots z_{j,M-k+1}%
\mathcal{D}_{i_{M-s}}.
\end{gather*}

\end{lemma}

\begin{proof}
We proceed by induction. The formula is tautological for $k=0$. Also the term
in the sum with $s=k-1$ has no $z_{jn}$ factors. Multiply the right hand side
by $z_{j,M-k}$ on the right. For the first part, $\mathcal{D}_{i_{M-k}}\left(
\left(  i_{M-k-1},i_{M-k}\right)  -C_{j,M-k}w_{M-k-1}\right)  =\left(
i_{M-k-1},i_{M-k}\right)  \mathcal{D}_{i_{M-k-1}}-C_{j,M-k}w_{M-k-1}%
\mathcal{D}_{i_{M-k}}$ (since $\mathcal{D}_{i_{n}}$ commutes with $w_{j}$ for
$n\neq j$). For the second part, $\mathcal{D}_{i_{M-s}}$ commutes with
$z_{j,M-k}$. This completes the induction.
\end{proof}

Set $k=M-j$ in the lemma, apply the operator to $\zeta_{\lambda}$ and multiply
both sides of the identity by $\left(  i_{j},i_{j+1}\right)  \ldots\left(
i_{M-1},i_{M}\right)  $,yielding (replace $s$ by $M-s$)
\begin{align*}
\mathcal{D}_{i_{j}}\zeta_{\lambda}  &  =\left(  i_{j},i_{j+1}\right)
\ldots\left(  i_{M-1},i_{M}\right)  \mathcal{D}_{m}\zeta_{\mu\left(
j,M\right)  }\\
&  +\sum_{s=j+1}^{M}C_{j,s}\left(  i_{j},i_{j+1}\right)  \ldots\left(
i_{s-1},i_{s}\right)  w_{s-1}z_{j,s-1}\ldots z_{j,j+1}\mathcal{D}_{i_{s}}%
\zeta_{\lambda}.
\end{align*}
This identity is used starting with $j=M-1$ and then decrementing $j$ by 1
with the result:
\[
\mathcal{D}_{i_{j}}\zeta_{\lambda}=\left(  i_{j},i_{j+1}\right)  \ldots\left(
i_{M-1},i_{M}\right)  \mathcal{D}_{m}\zeta_{\mu\left(  j,M\right)  }%
+\sum_{s=j+1}^{M}u_{j,s}\mathcal{D}_{m}\zeta_{\mu\left(  s,M\right)  },
\]
where each $u_{j,s}\in RS\left[  i_{1},m\right]  $ and $R$ is the $\mathbb{Z}%
$-ring generated by $\left\{  C_{j,k}:1\leq j<k\leq M\right\}  $. To complete
the analysis of $\mathcal{D}_{m}\zeta_{\mu\left(  j,M\right)  }$, for $0\leq
k<j$ let $\nu\left(  k,j\right)  \in\mathbb{N}_{0}^{N}$ be the action on
$\lambda-\varepsilon\left(  i_{j}\right)  $ by the (reverse) cyclic shift on
the interval $\left\{  i_{k}+1,\ldots,i_{j}\right\}  $, that is $\nu\left(
k,j\right)  _{i_{k}+1}=\lambda_{i_{j}}-1,\,\mu\left(  j,k\right)  _{i}%
=\lambda_{i-1}$ for $i_{k}+1<i\leq i_{j}$ and $\nu\left(  k,j\right)
_{i}=\lambda_{i}$ for $i\leq i_{k}$ or $i>i_{j}$. Also let $\nu\left(
j,j\right)  =\lambda-\varepsilon\left(  i_{j}\right)  $; if $i_{j}=i_{j-1}+1$
then $\nu\left(  j-1,j\right)  =\nu\left(  j,j\right)  $. For $0\leq k<j-1<M$
let
\begin{align*}
C_{kj}^{\prime}  &  =\frac{\kappa}{\left(  i_{j}-i_{k}-1\right)
\kappa+\lambda_{i_{k+1}}-\lambda_{i_{j}}+1},\\
w_{k}^{\prime}  &  =1+\sum_{r=i_{k}+2}^{i_{k+1}}\left(  r,i_{k+1}+1\right)
\in\mathbb{Z}S\left[  i_{k}+2,i_{k+1}+1\right]  ,
\end{align*}
and if $i_{j-1}<i_{j}-1$ let
\begin{align*}
C_{j-1,j}^{\prime}  &  =\frac{\kappa}{\left(  i_{j}-i_{j-1}-1\right)
\kappa+1},\\
w_{j-1}^{\prime}  &  =1+\sum_{r=i_{j-1}+2}^{i_{j}-1}\left(  r,i_{j}\right)
\in\mathbb{Z}S\left[  i_{j-1}+2,i_{j}\right]  .
\end{align*}
Proposition \ref{movelt} applies to the transformation of $\zeta_{\nu\left(
k,j\right)  }$ to $\zeta_{\nu\left(  k-1,j\right)  }$; note that\newline%
$r\left(  \nu\left(  k,j\right)  ,i_{k}+1\right)  =i_{k}+1$ and $r\left(
\nu\left(  k,j\right)  ,i_{k+1}\right)  =i_{j}$. Thus
\[
\zeta_{\nu\left(  j-1,j\right)  }=\left(  \left(  i_{j-1}+1,i_{j}\right)
-C_{j-1,j}^{\prime}w_{j-1}^{\prime}\right)  \zeta_{\nu\left(  j,j\right)  }%
\]
(unless $i_{j}=i_{j-1}+1$ when $\zeta_{\nu\left(  j-1,j\right)  }=\zeta
_{\nu\left(  j,j\right)  }$) and
\[
\zeta_{\nu\left(  k,j\right)  }=\left(  \left(  i_{k}+1,i_{k+1}+1\right)
-C_{kj}^{\prime}w_{k}^{\prime}\right)  \zeta_{\nu\left(  k+1,j\right)  }%
\]
for $0\leq k\leq j-2$. By Theorem \ref{difp} $\mathcal{D}_{m}\zeta_{\mu\left(
j,M\right)  }=\theta_{m}^{-1}\zeta_{\widetilde{\mu}\left(  j,M\right)  }$
where $\widetilde{\mu}\left(  j,M\right)  =\nu\left(  0,j\right)  $ (that is,
first the $i_{j}$-entry of $\lambda$ is moved to the $m$-entry at the end, the
action of $\mathcal{D}_{m}$ decrements $\lambda_{i_{j}}$ by 1 and moves it to
the front, loosely speaking). In turn $\zeta_{\widetilde{\mu}\left(
j,M\right)  }$ can be expressed in terms of $\zeta_{\nu\left(  j,j\right)  }$.
The following is now established.

\begin{theorem}
\label{bigdiff}Suppose $\lambda\in\mathbb{N}_{0}^{N,P}$ with points of
increase $i_{1}<i_{2}<\ldots<i_{M}=\ell\left(  \lambda\right)  $, let $R$ be
the $\mathbb{Z}$-ring generated by $\left\{  C_{jk}:1\leq j<k\leq M\right\}
\cup\allowbreak\left\{  C_{kj}^{\prime}:0\leq k<j\leq M\right\}
\cup\allowbreak\mathbb{Z}$ and let $\lambda^{\left(  j\right)  }%
=\lambda-\varepsilon\left(  i_{j}\right)  \in\mathbb{N}_{0}^{N,P}$, then for
$i_{j-1}<i\leq i_{j}$ with $0\leq j\leq M$,
\[
\mathcal{D}_{i}\zeta_{\lambda}=\sum_{s=j}^{M}\left(  \left(  N+1-i_{s}\right)
\kappa+\lambda_{i_{s}}\right)  u_{is}\zeta_{\lambda^{\left(  s\right)  }},
\]
where each $u_{is}\in RS\left[  1,m\right]  $.
\end{theorem}

The Theorem exhibits the poles in the differentiation formula for the
$p$-monic polynomials. To convert this for use with $x$-monic polynomials
multiply $\zeta_{\lambda^{\left(  s\right)  }}$ by $\left(  h\left(
\lambda,1\right)  h\left(  \lambda^{\left(  s\right)  },\kappa+1\right)
\right)  /\left(  h\left(  \lambda^{\left(  s\right)  },1\right)  h\left(
\lambda,\kappa+1\right)  \right)  $, then the identity holds for $\zeta$
replaced by $\zeta^{x}$. The details are not worked out since in general there
is no significant simplification. In the next section this calculation will be
carried out for the singular polynomials.

\section{Existence of singular polynomials}

In this section we will show for certain $\lambda\in\mathbb{N}_{0}^{N,P}$ and
singular values $\kappa_{0}$ that $\zeta_{\lambda}^{x}$ has no poles at
$\kappa=\kappa_{0}$ and that $\mathcal{D}_{i}\zeta_{\lambda}^{x}=0$ for $1\leq
i\leq N$. It turns out that for $m=\ell\left(  \lambda\right)  $ the last
coefficient in the formula of Theorem \ref{bigdiff} satisfies $\left(
N+1-m\right)  \kappa_{0}+\lambda_{m}=0$ and $\zeta_{\lambda-\varepsilon\left(
m\right)  }$ has no poles at $\kappa_{0}$ in general. For the terms of type
$\zeta_{\lambda^{\left(  s\right)  }}^{x}$ the denominator expression
$h\left(  \lambda^{\left(  s\right)  },\kappa+1\right)  $ has a zero at
$\kappa=\kappa_{0}$ but the pole $\left(  \kappa-\kappa_{0}\right)  $ does not
appear for the restriction to $\mathbb{R}^{N}$, and this is the key fact. We
start with the isotypes of two-part partitions ($\tau=\left(  \mu
,N-\mu\right)  $).

\begin{proposition}
Let $\frac{N}{2}\leq\mu<N,\,\gcd\left(  m,\mu+1\right)  <\frac{\mu+1}{N-\mu
},\,\lambda=\left(  m^{N-\mu}\right)  $ then $h\left(  \lambda,1\right)
,\allowbreak\,h\left(  \lambda,\kappa+1\right)  ,\,\allowbreak h\left(
\lambda-\varepsilon\left(  N-\mu\right)  ,1\right)  $ and $h\left(
\lambda-\varepsilon\left(  N-\mu\right)  ,\kappa+1\right)  $ are nonzero when
evaluated at $\kappa=-\frac{m}{\mu+1}$.
\end{proposition}

\begin{proof}
For $1\leq i\leq N-\mu,1\leq j\leq m$ we have $h\left(  \lambda,t;i,j\right)
=m-j+t+\left(  N-\mu-i\right)  \kappa$. For $t=1,\kappa+1$ the sets of values
are $\left\{  i\kappa+j:0\leq i\leq N-\mu-1,1\leq j\leq m\right\}
,\,\allowbreak\left\{  i\kappa+j:1\leq i\leq N-\mu,1\leq j\leq m\right\}  $
respectively. It suffices to show that the second set does not contain 0 for
$\kappa=-\frac{m}{\mu+1}$. Suppose $-im+j\left(  \mu+1\right)  =0$ for some
(nonzero) $i,j$ and let $d=\gcd\left(  m,\mu+1\right)  $, then $\frac{\mu
+1}{d}|\,i\leq N-\mu$ which implies $\frac{\mu+1}{N-\mu}\leq d$, contrary to
the hypothesis. For $\lambda-\varepsilon\left(  N-\mu\right)  $ only the
hook-lengths in the last row and column change; $h\left(  \lambda
-\varepsilon\left(  N-\mu\right)  ,t;i,m\right)  =t+\left(  N-\mu-i\right)
\kappa$ and $h\left(  \lambda-\varepsilon\left(  N-\mu\right)  ,t;N-\mu
,j\right)  =m-j+t$ for $1\leq i<N-\mu$ and $1\leq j<m$. These values have
already been shown to be nonzero for $\kappa=-\frac{m}{\mu+1}$.
\end{proof}

Next we handle the case of three or more parts, for the isotype $\tau
=$\newline$\left(  s\left(  \mu+1\right)  +\mu,\mu^{l},\rho\right)  $. The
following is the central hypothesis for this section.

\begin{definition}
\label{lambda}For $\mu,l\geq1,s\geq0,1\leq\rho\leq\mu$ and $\gcd\left(
m,\mu+1\right)  =1$ let
\[
\Lambda\left(  \mu,s,l,\rho,m\right)  =\left(  \left(  m\left(  s+l+1\right)
\right)  ^{\rho},\left(  m\left(  s+l\right)  \right)  ^{\mu},\ldots,\left(
m\left(  s+1\right)  \right)  ^{\mu}\right)  ,
\]
a partition of length $l\mu+\rho$ which is associated to the singular value
$\kappa_{0}=-\frac{m}{\mu+1}$ and the $S_{N}$-representation of isotype
$\left(  s\left(  \mu+1\right)  +\mu,\mu^{l},\rho\right)  $, where $N=\left(
s+l+1\right)  \mu+s+\rho$.
\end{definition}

\begin{lemma}
\label{lmnon0}Suppose $a,b,c\in\mathbb{N}_{0}$ and $c\geq1,b\leq\mu$ then
$a\left(  \mu\kappa+m\right)  +b\kappa+c\neq0$ for $\kappa=-\frac{m}{\mu+1}$
(where $\gcd\left(  m,\mu+1\right)  =1$).
\end{lemma}

\begin{proof}
Denote the value of the expression at $\kappa=-\frac{m}{\mu+1}$by $v$, then
$v=\left(  a-b\right)  \frac{m}{\mu+1}+c$. If $a\geq b$ then $v\geq c\geq1$;
otherwise $0>a-b\geq-\mu$ and $\mu+1$ does not divide $\left(  a-b\right)  m$
thus $v\notin\mathbb{Z}$ and $v\neq0$.
\end{proof}

\begin{proposition}
Let $\lambda=\Lambda\left(  \mu,s,l,\rho,m\right)  $, then $h\left(
\lambda,1\right)  $ and $h\left(  \lambda,\kappa+1\right)  $ are nonzero when
$\kappa=-\frac{m}{\mu+1}$.
\end{proposition}

\begin{proof}
Since hook-lengths in a given row depend only on it and the rows of higher
index we may assume that $\rho=\mu$. We index the rows of $\lambda$ by $\mu
k+i$ with $0\leq k\leq l$ and $1\leq i\leq\mu$, and the columns by $m\left(
s+n\right)  -j$ where $1\leq n\leq l+1-k$ and $0\leq j\leq m-1$ if $n>1$, or
$0\leq j<m\left(  s+1\right)  $ if $n=1$. Then
\begin{align*}
h\left(  \lambda,t;\mu k+i,m\left(  s+n\right)  -j\right)   &  =\kappa\left(
\left(  l+2-k-n\right)  \mu-i\right)  +m\left(  l+1-k-n\right)  +j+t\\
&  =\left(  l+1-k-n\right)  \left(  \kappa\mu+m\right)  +\kappa\left(
\mu-i\right)  +j+t.
\end{align*}
Set $i^{\prime}=\mu-i+1$ if $t=\kappa+1$ or $i^{\prime}=\mu-i$ if $t=1$; then
the above expression equals $\left(  l+1-k-n\right)  \left(  \kappa
\mu+m\right)  +\kappa i^{\prime}+j+1$ which is nonzero at $\kappa=\kappa_{0}$
by Lemma \ref{lmnon0} (since $i^{\prime}\leq\mu$).
\end{proof}

Next we consider the hook-lengths for $\Lambda\left(  \mu,s,l,\rho,m\right)
-\varepsilon\left(  \rho+k\mu\right)  $ with $0\leq k\leq l$.

\begin{proposition}
For $0\leq k_{0}\leq l$ let $\nu=\Lambda\left(  \mu,s,l,\rho,m\right)
-\varepsilon\left(  \rho+k_{0}\mu\right)  $, then for $\kappa=-\frac{m}{\mu
+1}$ $h\left(  \nu,1\right)  $ is nonzero and $h\left(  \nu,\kappa+1\right)  $
is nonzero for $k_{0}=l$ and has a zero of multiplicity one for $0\leq
k_{0}<l$, in the hook-length $h\left(  \nu,\kappa+1;\rho+k_{0}\mu,m\left(
s+l+1-k_{0}\right)  \right)  $.
\end{proposition}

\begin{proof}
As in the previous proof, assume $\rho=\mu$. The column above the node deleted
from $\lambda$ (namely, $m\left(  s+l+1-k_{0}\right)  $) meets the rows
labeled by $\mu k+i$ with $0\leq k\leq k_{0}$ and $1\leq i\leq\mu$, except
$1\leq i<\mu$ when $k=k_{0}$. Then $h\left(  \nu,t;\mu k+i,m\left(
s+l+1-k_{0}\right)  \right)  =\left(  k_{0}-k\right)  \left(  \kappa
\mu+m\right)  +\kappa i^{\prime}+1$ where $i^{\prime}=\mu-1-i$ for $t=1$ and
$i^{\prime}=\mu-i$ for $t=\kappa+1$. By Lemma \ref{lmnon0} the value is
nonzero for $\kappa=\kappa_{0}$. The row of the deleted node meets the columns
labeled $m\left(  s+n\right)  -j$ with $1\leq n\leq l+1-k_{0}$. Then $h\left(
\nu,t;m\left(  s+n\right)  -j,\left(  k_{0}+1\right)  \mu\right)  =\left(
l+1-k_{0}-n\right)  \left(  \kappa\mu+m\right)  +b\kappa+j$, where $b=0$ for
$t=1$ and $b=1$ for $t=\kappa+1$. The Lemma applies unless $j=0$. Suppose
$j=0$ then $1\leq n\leq l-k_{0}$ (the value $j=0$ does not occur for
$n=l+1-k_{0}$ since the corresponding node was deleted); at $\kappa=\kappa
_{0}$ the value of the hook-length is $\left(  l+1-k_{0}-n-b\right)  \frac
{m}{\mu+1}$ which is zero exactly when $n=l-k_{0}$ and $b=1$ (that is,
$t=\kappa+1$). Thus the hook-length $h\left(  \nu,\kappa+1;m\left(
s+l-k_{0}\right)  ,\left(  k_{0}+1\right)  \mu\right)  =\kappa\left(
\mu+1\right)  +m$ is the only zero in $h\left(  \nu,\kappa+1\right)  $ at
$\kappa=\kappa_{0}$.
\end{proof}

Next we show that the coefficients $C_{jk}$ and $C_{kj}^{\prime}$ appearing in
Theorem \ref{bigdiff} have no poles at $\kappa=\kappa_{0}$. The points of
decrease of $\Lambda\left(  \mu,s,l,\rho,m\right)  $ are $i_{j}=\rho+\left(
j-1\right)  \mu$, $\lambda_{i_{j}}=m\left(  s+2+l-j\right)  $ for $1\leq j\leq
l+1$. For $1\leq j<k\leq l+1$ the coefficient $C_{jk}=\frac{\kappa}{\left(
i_{k}-i_{j}\right)  \kappa+\lambda_{i_{j}}-\lambda_{i_{k}}}=\frac{\kappa
}{\left(  k-j\right)  \left(  \kappa\mu+m\right)  }$, which has value
$-\frac{1}{k-j}$ at $\kappa=\kappa_{0}$.

\begin{proposition}
For $\lambda=\Lambda\left(  \mu,s,l,\rho,m\right)  $ the coefficients
$C_{kj}^{\prime}$ have no poles at $\kappa=-\frac{m}{\mu+1}$ for $0\leq
k<j\leq l+1$.
\end{proposition}

\begin{proof}
First the special cases $C_{j-1,j}^{\prime}=\frac{\kappa}{\left(
\mu-1\right)  \kappa+1}$ for $j>0,$\thinspace$\mu>1$ and $C_{0,1}^{\prime
}=\frac{\kappa}{\left(  \rho-1\right)  \kappa+1}$ for $\rho>1$ are obviously
finite at $\kappa=\kappa_{0}$. Next for $j>1$ we have $C_{0,j}^{\prime}%
=\frac{\kappa}{\left(  \left(  j-1\right)  \mu+\rho-1\right)  \kappa+m\left(
j-1\right)  +1}$ with denominator $\left(  j-1\right)  \left(  \mu
\kappa+m\right)  +\left(  \rho-1\right)  \kappa+1$ which is nonzero at
$\kappa=\kappa_{0}$ by Lemma \ref{lmnon0}$,$ since $\rho-1\leq\mu$. Finally
for $1\leq k<j-1\leq l$ we have $C_{k,j}^{\prime}=\kappa\left(  \left(
j-k-1\right)  \left(  \mu\kappa+m\right)  +\left(  \mu-1\right)
\kappa+1\right)  ^{-1}$, and the Lemma applies.
\end{proof}

We restate the result of Theorem \ref{bigdiff} applied to the $x$-monic
polynomials $\zeta_{\lambda}^{x}$ and $\zeta_{\lambda-\varepsilon\left(
\rho+k\mu\right)  }^{x}$ (for $\lambda=\Lambda\left(  \mu,s,l,\rho,m\right)  $
and $0\leq k\leq l$). For $1\leq i\leq\rho+l\mu$
\begin{align*}
\mathcal{D}_{i}\zeta_{\lambda}^{x}  &  =\sum_{k=0}^{l}\left(  \left(
l-k\right)  \left(  \kappa\mu+m\right)  +\left(  s+1\right)  \left(
m+\kappa\left(  \mu+1\right)  \right)  \right)  u_{i,k+1}\\
&  \times\frac{h\left(  \lambda,1\right)  h\left(  \lambda-\varepsilon\left(
\rho+k\mu\right)  ,\kappa+1\right)  }{h\left(  \lambda,\kappa+1\right)
h\left(  \lambda-\varepsilon\left(  \rho+k\mu\right)  ,1\right)  }%
\zeta_{\lambda-\varepsilon\left(  \rho+k\mu\right)  }^{x},
\end{align*}
where (the labeling of the points of decrease is now shifted by 1) each
$u_{i,k+1}\in RS\left[  1,\rho+l\mu\right]  $ and $R$ is the ring generated by
$\left\{  C_{jk}:1\leq j<k\leq l+1\right\}  \cup$\newline$\left\{
C_{kj}^{\prime}:0\leq k<j\leq l+1\right\}  \cup\mathbb{Z}$; also $u_{i,k+1}=0$
for $k<\frac{i-\rho}{\mu}$. Since $h\left(  \lambda,\kappa+1\right)  \neq0$ at
$\kappa=\kappa_{0}$ the polynomial $\zeta_{\lambda}$ has no poles there. Also
the specialization of $R$ is a subring of $\mathbb{Q}$. For $k=l$ we already
have shown that $h\left(  \lambda-\varepsilon\left(  \rho+l\mu\right)
,\kappa+1\right)  \neq0$ and thus $\zeta_{\lambda-\varepsilon\left(  \rho
+l\mu\right)  }^{x}$ has no poles at $\kappa_{0}$ and the factor $\left(
s+1\right)  \left(  m+\kappa\left(  \mu+1\right)  \right)  $ becomes zero.
When $0\leq k<l$ the factor $h\left(  \lambda-\varepsilon\left(  \rho
+k\mu\right)  ,\kappa+1\right)  $ has a zero at $\kappa_{0}$. Once we prove
that $\zeta_{\lambda-\varepsilon\left(  \rho+k\mu\right)  }^{x}$ has no pole
at $\kappa_{0}$ the proof that $\mathcal{D}_{i}\zeta_{\lambda}=0$ for all $i$
will be complete.

The method of Knop and Sahi \cite{KS} was designed to show that the
coefficients of the monomials $x^{\beta}$ in $h\left(  \lambda,\kappa
+1\right)  \zeta_{\lambda}^{x}$ are in $\mathbb{N}_{0}\left[  \kappa\right]
$, but it is not evident how to use the method to identify the poles when the
number of variables is in the range $\ell\left(  \lambda\right)  \leq N<$
$\ell\left(  \lambda\right)  +\left\vert \lambda\right\vert $. We introduce a
different approach.

\begin{definition}
\label{critpair}Let $\alpha,\beta\in\mathbb{N}_{0}^{M}$ with $\alpha
\vartriangleright\beta$ and let $m,n\in\mathbb{N}$ with $\gcd\left(
m,n\right)  =1$ then say $\left(  \alpha,\beta\right)  $ is a $\left(
-\frac{m}{n}\right)  $-critical pair if $\left(  n\kappa+m\right)  $ divides
$\left(  r\left(  \beta,i\right)  -r\left(  \alpha,i\right)  \right)
\kappa+\alpha_{i}-\beta_{i}$ (in $\mathbb{Q}\left[  \kappa\right]  $) for
$1\leq i\leq M$.
\end{definition}

In fact the division is in $\mathbb{Z}\left[  \kappa\right]  $ because
$\gcd\left(  m,n\right)  =1$. The definition will be used in the situation
$\alpha\in\mathbb{N}_{0}^{N,P}$ that is, $\ell\left(  \alpha\right)  \leq N$
and $M=\ell\left(  \alpha\right)  +\left|  \alpha\right|  $. See Definitions
\ref{order} and \ref{rankdef} for the order $\vartriangleright$ and the rank
function $r$.

\begin{theorem}
Suppose $\lambda\in\mathbb{N}_{0}^{N,P}$ and $\kappa_{0}\in\mathbb{Q}%
,\kappa_{0}<0$; if there does not exist $\beta\in\mathbb{N}_{0}^{N}$ such that
$\left(  \lambda,\beta\right)  $ is a $\kappa_{0}$-critical pair then
$\kappa_{0}$ is not a pole of $\zeta_{\lambda}^{x}$ restricted to
$\mathbb{R}^{N}$.
\end{theorem}

\begin{proof}
Extend the field $\mathbb{Q}\left(  \kappa\right)  $ with $N$ transcendental
variables $\left\{  v_{1},v_{2},\ldots,v_{N}\right\}  $ and let $\mathcal{T}%
=\sum_{i=1}^{N}v_{i}\mathcal{U}_{i}.$ For each $\alpha\in\mathbb{N}_{0}^{N} $
the polynomial $\zeta_{\alpha}^{x}$ is an eigenvector of $\mathcal{T}$, indeed
$\mathcal{T}\zeta_{\alpha}^{x}=\sum_{i=1}^{N}v_{i}\xi_{i}\left(
\alpha\right)  \zeta_{\alpha}^{x}$. The eigenvalue determines $\alpha$
uniquely for generic $\kappa$ (with the possible exception of a finite set of
negative rationals). Let $C=\left\{  \beta\in\mathbb{N}_{0}^{N}:\lambda
\vartriangleright\beta\right\}  $. By the triangularity of the operators
$\left\{  \mathcal{U}_{i}\right\}  $ we have $x^{\lambda}=\zeta_{\lambda}%
^{x}+\sum_{\beta\in C}B_{\beta\lambda}\zeta_{\beta}^{x}$ for certain
coefficients $B_{\beta\lambda}\in\mathbb{Q}\left(  \kappa\right)  $. Let
\[
\mathcal{T}_{\lambda}=\prod_{\beta\in C}\frac{\mathcal{T}-\sum_{i=1}^{N}%
v_{i}\xi_{i}\left(  \beta\right)  }{\sum_{i=1}^{N}v_{i}\left(  \xi_{i}\left(
\lambda\right)  -\xi_{i}\left(  \beta\right)  \right)  },
\]
then $\mathcal{T}_{\lambda}x^{\lambda}=\zeta_{\lambda}^{x}$ (note that the
number $N$ of variables is part of the definition of the set $C$). The
numerator of the product is a polynomial in $\kappa,v_{1},\ldots,v_{N}$ (and
of course each $\mathcal{D}_{i}x^{\alpha}$ is a polynomial with coefficients
in $\mathbb{Z}\left[  \kappa\right]  $) thus any ($\kappa$)-poles in
$\zeta_{\lambda}^{x}$ must appear in the set $\left\{  \sum_{i=1}^{N}%
v_{i}\left(  \xi_{i}\left(  \lambda\right)  -\xi_{i}\left(  \beta\right)
\right)  :\lambda\vartriangleright\beta\right\}  $. For any $\beta\in C$ we
have $\sum_{i=1}^{N}v_{i}\left(  \xi_{i}\left(  \lambda\right)  -\xi
_{i}\left(  \beta\right)  \right)  =\sum_{i=1}^{N}v_{i}\left(  \left(
r\left(  \beta,i\right)  -r\left(  \lambda,i\right)  \right)  \kappa
+\lambda_{i}-\beta_{i}\right)  $. Since any denominator appearing in a
coefficient (with respect to the $x$-monomial basis) of $\zeta_{\lambda}^{x}$
must be a factor of $h\left(  \lambda,\kappa+1\right)  $, all of the terms
involving $\left\{  v_{i}:1\leq i\leq N\right\}  $ must cancel out in the
calculation of $\mathcal{T}_{\lambda}x^{\lambda}$. Thus the irreducible
polynomials $\sum_{i=1}^{N}v_{i}\left(  \xi_{i}\left(  \lambda\right)
-\xi_{i}\left(  \beta\right)  \right)  $ must cancel out and the denominators
in $\zeta_{\lambda}^{x}$ can only arise from reducible terms of the form
$\left(  \sum_{i=1}^{N}a_{i}v_{i}\right)  \left(  \kappa-\kappa_{1}\right)  $
where $a_{1},\ldots,a_{N},\kappa_{1}\in\mathbb{Q}$. This condition is
equivalent to $\left(  \lambda,\beta\right)  $ being a $\kappa_{1}$-critical
pair. Thus, if there is no $\kappa_{0}$-critical pair $\left(  \lambda
,\beta\right)  $ with $\ell\left(  \beta\right)  \leq N$ then $\kappa_{0}$ is
not a pole of $\zeta_{\lambda}^{x}.$
\end{proof}

We will exploit this theorem by directly constructing the unique $\beta$ such
that $\left(  \Lambda\left(  \mu,s,l,\rho,m\right)  -\varepsilon\left(
\rho+k\mu\right)  ,\beta\right)  $ is $-\frac{m}{\mu+1}$-critical. Here is a
numerical example: for $N=33$ consider $\Lambda\left(  3,4,4,2,3\right)  $ for
the singular value $\kappa_{0}=-\frac{3}{4}$ of isotype $\left(
19,3^{4},2\right)  $, take $k=1$, then $\lambda=\left(  27^{2},24^{2}%
,23,21^{3},18^{3},15^{3}\right)  $ and the unique $\beta$ such that $\left(
\lambda,\beta\right)  $ is $\left(  -\frac{3}{4}\right)  $-critical is
$\left(  27^{2},24^{2},2,0^{3},21^{3},18^{3},3^{22}\right)  $. The
construction proceeds through several lemmas. Fix $k$ such that $0\leq k\leq
l-1$, let $\lambda=\Lambda\left(  \mu,s,l,\rho,m\right)  -\varepsilon\left(
\rho+k\mu\right)  ,L=\ell\left(  \lambda\right)  =\rho+l\mu$; and partition
$\left[  1,L\right]  $ into cells $\left\{  I_{j}:0\leq j\leq l\right\}  $,
where $I_{0}=\left[  1,\rho\right]  $ and $I_{j}=\left[  \rho+\left(
j-1\right)  \mu+1,\rho+j\mu\right]  $ for $1\leq j\leq l$. Then $i\in I_{j}$
implies $\lambda_{i}=m\left(  l+s+1-j\right)  $, except that $\lambda
_{\rho+k\mu}=m\left(  l+s+1-k\right)  -1$. We will show the required $\beta$
has the values $m\left(  l+s+1-j\right)  $ on the cells $I_{j}$ with $j\leq
k$, $0$ on $I_{k+1}$, $m\left(  l+s+2-j\right)  $ on $I_{j}$ with $k<j\leq l$,
except $\beta_{\rho+k\mu}=m-1$, and $\beta_{i}=m$ for $L+1\leq i\leq
L+\mu\left(  s+1\right)  +l+s-k=N+l-k$; also $\ell\left(  \beta\right)
=N+l-k$. Henceforth, suppose that $\left(  \lambda,\beta\right)  $ is
$-\frac{m}{\mu+1}$-critical or $\beta=\lambda$. This holds if and only if the
\textit{rank equation}
\begin{equation}
r\left(  \beta,i\right)  -i=\left(  \mu+1\right)  \left(  \frac{1}{m}\left(
\lambda_{i}-\beta_{i}\right)  \right)  \label{keyq}%
\end{equation}
is satisfied for all $i\geq1$. Since $\gcd\left(  m,\mu+1\right)  =1$ this
implies that $\beta_{i}\equiv\lambda_{i}\operatorname{mod}m$ (so with the
exception of $\beta_{\rho+k\mu}$ each $\beta_{i}$ is divisible by $m$). Here
is a maximum principle for the multiplicity $\#\left\{  j:\beta_{j}%
=m\gamma,1\leq j\leq L\right\}  $ for any $\gamma$. There is a slight
difference for the cases $m=1$ and $m>1$. The condition $\lambda
\trianglerighteq\beta$ implies that any possible values satisfy $\gamma\leq
s+l+1$.

\begin{lemma}
\label{maxppm}Suppose $\gamma\in\mathbb{N}_{0},$ and $G=\left\{  j:\beta
_{j}=m\gamma,1\leq j\leq L\right\}  $, $m>1$ or $m=1$ and $\rho+k\mu\notin G$,
if $G$ meets two or more cells then $\#G\leq\mu-1$; additionally, if one of
the cells is $I_{0}$ then $\#G\leq\rho-1$.
\end{lemma}

\begin{proof}
Let $G$ have nonempty intersections with cells $I_{g_{1}},I_{g_{2}}%
,\ldots,I_{g_{u}}$ with $0\leq g_{1}<g_{2}<\ldots<g_{u}\leq l$. By hypothesis
$\rho+k\mu\notin G$ (if $m>1$ then $m$ does not divide $\beta_{\rho+k\mu} $)
and so $i\in G\cap I_{g_{a}}$ implies $\lambda_{i}=m\left(  s+l+1-g_{a}%
\right)  $. Each $G\cap I_{g_{a}}$ is an interval $\left[  i_{a},j_{a}\right]
$; indeed suppose $i,j\in G\cap I_{g_{a}}$ and $i<j$, then by equation
(\ref{keyq}) $r\left(  \beta,i\right)  -i=\left(  \mu+1\right)  \left(
\frac{1}{m}\left(  m\left(  s+l+1-g_{a}\right)  -m\gamma\right)  \right)
=r\left(  \beta,j\right)  -j$; thus $r\left(  \beta,j\right)  =r\left(
\beta,i\right)  +j-i$. Since $\beta_{j}=\beta_{i}$ this implies that
$\beta_{b}=\beta_{i}=m\gamma$ for $i\leq b\leq j$. For $0\leq a<u$ we have
that $r\left(  \beta,i_{a+1}\right)  =r\left(  \beta,j_{a}\right)  +1$ and we
combine the two equations
\begin{align*}
r\left(  \beta,i_{a+1}\right)  -i_{a+1}  &  =\left(  \mu+1\right)  \left(
l+s+1-g_{a+1}-\gamma\right)  ,\\
r\left(  \beta,j_{a}\right)  -j_{a}  &  =\left(  \mu+1\right)  \left(
l+s+1-g_{a}-\gamma\right)
\end{align*}
to obtain $i_{a+1}=j_{a}+1+\left(  \mu+1\right)  \left(  g_{a+1}-g_{a}\right)
$. Then $\#G=\sum_{a=1}^{u}\left(  j_{a}-i_{a}+1\right)  =u+j_{u}-i_{1}%
-\sum_{a=1}^{u-1}\left(  i_{a+1}-j_{a}\right)  =\allowbreak u+j_{u}-i_{1}%
-\sum_{a=1}^{u-1}\left(  1+\left(  \mu+1\right)  \left(  g_{a+1}-g_{a}\right)
\right)  =1+j_{u}-i_{1}-\left(  \mu+1\right)  \left(  g_{u}-g_{1}\right)  $.
But $j_{u}\leq\rho+g_{u}\mu$ and $i_{1}\geq\rho+\left(  g_{1}-1\right)  \mu+1$
for $g_{1}\geq1$ while $i_{1}\geq1$ for $g_{1}=0$. This shows that $\#G\leq
\mu-\left(  g_{u}-g_{1}\right)  $ for $g_{1}\geq1$ and $\#G\leq\rho-g_{u}$ if
$g_{1}=0$. In both cases $\#G\leq\mu-1$.
\end{proof}

\begin{lemma}
\label{maxpp1}Suppose $m=1$,$\,\gamma\in\mathbb{N}_{0},$ $G=\left\{
j:\beta_{j}=\gamma,1\leq j\leq L\right\}  $, and $G\backslash\left\{
\rho+k\mu\right\}  $ meets two or more cells, if $\rho+k\mu=\min\left(
G\right)  $ and $G\neq\left[  \rho+k\mu,\rho+\left(  k+1\right)  \mu\right]  $
then $\#G\leq\mu$, otherwise$\ (\rho+k\mu\neq\min\left(  G\right)  $) then
$\#G\leq\mu-1$.
\end{lemma}

\begin{proof}
By hypothesis $G\neq\left[  \rho+k\mu,\rho+\left(  k+1\right)  \mu\right]  $.
We can apply the previous argument if we replace $I_{k}$ by $I_{k}%
\backslash\left\{  \rho+k\mu\right\}  $ and $I_{k+1}$ by $I_{k+1}\cup\left\{
\rho+k\mu\right\}  =$\newline$\left[  \rho+k\mu,\rho+\left(  k+1\right)
\mu\right]  $. If $g_{1}\neq k+1$ then as before $\#G\leq\mu-\left(
g_{u}-g_{1}\right)  \leq\mu-1$. If $g_{1}=k+1$ and $i_{1}\geq\rho+k\mu+1$ the
same conclusion results. When $g_{1}=k+1$ and $i_{1}=\rho+k\mu$, that is,
$\min G=\rho+k\mu$, the calculation yields the bound $\#G\leq\mu+1-\left(
g_{u}-k-1\right)  \leq\mu$.
\end{proof}

The two lemmas show that $\#\left\{  j:\beta_{j}=\gamma m,1\leq j\leq
L\right\}  \leq\mu$ for any $\gamma\in\mathbb{N}_{0}$, except when $m=1$ and
$\left\{  j:\beta_{j}=\gamma\right\}  =\left[  \rho+k\mu,\rho+\left(
k+1\right)  \mu\right]  $ of cardinality $\mu+1$. Next we show $\beta
_{L+1}\leq m$.

\begin{lemma}
Either $\beta_{L+1}=m$ and $r\left(  \beta,L+1\right)  =L-\mu$, or
$\beta_{L+1}=0$ and $r\left(  \beta,L+1\right)  =L+1$, $\ell\left(
\beta\right)  =L $.
\end{lemma}

\begin{proof}
Denote $\frac{\beta_{L+1}}{m}$ by $b$; note that $b\in\mathbb{N}_{0}$.First we
show $b\leq l$: by equation \ref{keyq} $r\left(  \beta,L+1\right)
=L+1-\left(  \mu+1\right)  b\geq1$ and so $b\leq\frac{\rho+l\mu}{\mu+1}%
\leq\left(  l+1\right)  \frac{\mu}{\mu+1}<l+1$. Let $a_{0}=\#\left\{  j:1\leq
j\leq L,\beta_{j}<\beta_{L+1}\right\}  $ and $a_{1}=\#\left\{  j:j>L,\beta
_{j}>\beta_{L+1}\right\}  $ then $r\left(  \beta,L+1\right)  =L+1-a_{0}%
+a_{1}\geq L+1-a_{0}.$ We claim $a_{0}\leq b\mu+1$. If $m>1$ then $a_{0}%
=\sum_{i=0}^{b-1}\#\left\{  j:\beta_{j}=im,1\leq j\leq L\right\}
+\allowbreak\#\left\{  j:\beta_{j}=cm-1,c\leq b\right\}  $. By the maximum
principle $a_{0}\leq b\mu+1$. If $m=1$ then $a_{0}=\sum_{i=0}^{b-1}\#\left\{
j:\beta_{j}=i,1\leq j\leq L\right\}  $; at most one of these sets can have
cardinality $\mu+1$ and again $a_{0}\leq b\mu+1$. Then $L+1-\left(
\mu+1\right)  b=r\left(  \beta,L+1\right)  \geq L-b\mu$ , that is, $b\leq1$.
If $b=1$ then $r\left(  \beta,L+1\right)  =L-\mu$. If $b=0$ then $r\left(
\beta,L+1\right)  =L+1$ which implies $\beta_{j}=0$ for all $j>L$. The
hypothesis $\lambda\trianglerighteq\beta$ implies $L=\ell\left(
\lambda\right)  \leq\ell\left(  \beta\right)  $.
\end{proof}

In fact, $\beta_{L+1}=0$ implies $\beta=\lambda$ and $\beta_{L+1}=m$
corresponds to a unique solution $\beta$ with $\ell\left(  \beta\right)
=N+l-k$.

\begin{lemma}
Suppose that $\beta_{L+1}=m$ then $\beta_{\rho+k\mu}=m-1$, $\beta_{i}=0$ for
$i\in I_{k+1}$, $\beta_{i}=m$ for $L+1\leq i\leq N+l-k$ and $\ell\left(
\beta\right)  =N+l-k$.
\end{lemma}

\begin{proof}
Let $a_{0}=\#\left\{  j:1\leq j\leq L,\beta_{j}\geq m\right\}  ,\,a_{1}%
=\#\left\{  j:L<j,\beta_{j}>m\right\}  ,$\newline$\,G_{0}=\left\{  j:1\leq
j\leq L,\beta_{j}=0\right\}  ,\allowbreak\,G_{1}=\left\{  j:\beta
_{j}=m-1>0\right\}  ,$ and $a_{2}=\#G_{0}+\#G_{1}$ where $G_{1}$ is empty when
$m=1$. Then $L-\mu=r\left(  \beta,L+1\right)  =a_{0}+a_{1}+1$ and
$L=a_{0}+a_{2}$, thus $a_{2}=\mu+1+a_{1}\geq\mu+1$. But by the maximum
principle $a_{2}\leq\mu+1$, hence $a_{1}=0$ and $a_{2}=\mu+1$. If $m>1$ then
$\#G_{0}\leq\mu$ implying that $\#G_{1}=1$ and $G_{1}=\left\{  \rho
+k\mu\right\}  $, also $\#G_{0}=\mu$ and thus $G_{0}=I_{j}$ for some
$j\neq0,k$ by Lemma \ref{maxppm}. If $m=1$ then $G_{0}=\left[  \rho+k\mu
,\rho+\left(  k+1\right)  \mu\right]  $ by Lemma \ref{maxpp1}. Let
$r_{0}=r\left(  \beta,\rho+k\mu\right)  $, then $r\left(  \beta,\rho+\left(
j-1\right)  \mu+1\right)  =r_{0}+1$; if $m=1$ then $\beta_{\rho+k\mu}%
=\beta_{\rho+k\mu+1}=0$ and $j=k+1$, while for $m>1$ we have $\beta_{\rho
+k\mu}=m-1$ which is the unique minimum of $\left\{  \beta_{j}:1\leq j\leq
L,\beta_{j}>0\right\}  $ and $\beta_{\rho+\left(  j-1\right)  \mu+1}$ is the
first occurrence of 0. By equation \ref{keyq}
\begin{align*}
r_{0}-\left(  \rho+k\mu\right)   &  =\frac{\mu+1}{m}\left(  m\left(
s+l+1-k\right)  -1-\left(  m-1\right)  \right) \\
&  =\left(  \mu+1\right)  \left(  s+l-k\right)  ,\\
r_{0}+1-\left(  \rho+\left(  j-1\right)  \mu+1\right)   &  =\left(
\mu+1\right)  \left(  s+l+1-j\right)  .
\end{align*}
Thus $r_{0}=\rho+k\mu+\left(  \mu+1\right)  \left(  s+l-k\right)  $ and
$\left(  j-1-k\right)  \mu=\left(  \mu+1\right)  \left(  j-1-k\right)  $, that
is $j=k+1$. But $r_{0}=\#\left\{  j:\beta_{j}\geq m\right\}  +1=L-\mu
+\#\left\{  j:L<j,\beta_{j}=m\right\}  $ and so $\#\left\{  j:L<j,\beta
_{j}=m\right\}  =\mu s+\mu+s+l-k=N-L+l-k$. This shows $\ell\left(
\beta\right)  =N+l-k\geq N+1$.
\end{proof}

Certainly this, together with a proof that $\ell\left(  \beta\right)  =L$
implies $\beta=\lambda,$ is enough for the main purpose, but with not much
more work we can show that $\beta$ is unique. In fact we will show that
$\beta_{L+1}=m$ implies that for $i\in I_{j}$ $\beta_{i}=m\left(
s+l+1-j\right)  $ for $j\leq k$ and $\beta_{i}=m\left(  s+l+2-j\right)  $ for
$j>k+1$, except $\beta_{\rho+k\mu}=m-1$.

\begin{lemma}
Suppose that $\beta_{L+1}=m$, then $\beta_{i}=\lambda_{i}$ for all
$i<\rho+k\mu$ and $\beta_{i}=\lambda_{i}+m$ for $\rho+(k+1)\mu+1\leq i\leq
\ell\left(  \beta\right)  $.
\end{lemma}

\begin{proof}
For $0\leq i\leq s+l$ let $M_{i}=\#\left\{  j:1\leq j\leq L,\beta_{j}=m\left(
s+l+1-i\right)  \right\}  $. Since $\lambda_{i}=m\left(  s+l+1\right)  $ for
$1\leq i\leq\rho$ when $k\geq1$, and for $1\leq i\leq\rho-1$ when $k=0$, the
condition $\lambda\vartriangleright\beta$ (thus $\lambda\succeq\beta^{+}$)
implies $M_{0}\leq\rho$ or $\rho-1$ respectively. The maximum principle
(Lemmas \ref{maxppm} and \ref{maxpp1}) implies that $M_{i}\leq\mu$ for $1\leq
i<l+1$ (from the previous lemma in which $\beta_{\rho+k\mu}=m-1$ was
determined). Further $\sum_{i=0}^{s+l}M_{i}=L-\mu-1=\rho+\left(  l-1\right)
\mu-1$, that is, $\sum_{i=l}^{s+l}M_{i}=\rho-M_{0}+\sum_{i=1}^{l-1}\left(
\mu-M_{i}\right)  -1$. Also
\begin{align*}
\left\vert \beta\right\vert  &  =\sum_{i=0}^{s+l}M_{i}m\left(  s+l+1-i\right)
+m-1+m\left(  \mu s+\mu+s+l-k\right) \\
&  =\rho m\left(  s+l+1\right)  +m\mu\sum_{i=1}^{l}\left(  s+i\right)  -1,
\end{align*}
and so
\begin{align*}
\sum_{i=l}^{s+l}M_{i}\left(  s+l+1-i\right)   &  =\left(  \rho-M_{0}\right)
\left(  s+l+1\right) \\
&  +\sum_{i=1}^{l-1}\left(  \mu-M_{i}\right)  \left(  s+l+1-i\right)  -\left(
s+l+1-k\right)  .
\end{align*}
Let $j$ be defined by $M_{0}=\rho,M_{i}=\mu$ for $1\leq i\leq j$ and
$M_{j+1}\leq\mu-1$, that is, $j\geq0$, while $j=-1$ when $M_{0}\leq\rho-1$.
The hypothesis $\lambda\succeq\beta^{+}$ implies $j<k$ (or else $\sum
_{i=1}^{\rho+k\mu}\beta_{i}^{+}>\sum_{i=1}^{\rho+k\mu}\lambda_{i}$). Then for
$j=-1$ we have
\begin{align*}
\sum_{i=l}^{s+l}M_{i}\left(  s+l+1-i\right)   &  =\left(  \rho-1-M_{0}\right)
\left(  s+l+1\right)  +\sum_{i=1}^{l-1}\left(  \mu-M_{i}\right)  \left(
s+l+1-i\right)  +k\\
&  \geq\left(  s+2\right)  \left(  \rho-1-M_{0}+\sum_{i=1}^{l-1}\left(
\mu-M_{i}\right)  \right)  =\sum_{i=l}^{s+l}M_{i}\left(  s+2\right)  ,
\end{align*}
and for $j\geq0$%
\begin{align*}
\sum_{i=l}^{s+l}M_{i}\left(  s+l+1-i\right)   &  =\left(  \rho-M_{0}\right)
\left(  s+l+1\right)  +\sum_{i=1,i\neq j+1}^{l-1}\left(  \mu-M_{i}\right)
\left(  s+l+1-i\right) \\
&  +\left(  \mu-1-M_{j+1}\right)  \left(  s+l-j\right)  +\left(  k-1-j\right)
\\
&  \geq\sum_{i=l}^{s+l}M_{i}\left(  s+2\right)  ,
\end{align*}
since each coefficient is nonnegative, but then $M_{i}=0$ for all $i\geq l$.
The nonnegativity of each term on the right hand sides implies $j=k-1$ and if
$k=0$ then $M_{0}=\rho-1$ and $M_{i}=\mu$ for $1\leq i\leq l-1$, or else
$M_{0}=\rho,M_{k}=\mu-1$ and $M_{i}=\mu$ for $1\leq i\leq l-1,i\neq k$. Let
$G_{i}=\left\{  j:1\leq j\leq L,\beta_{j}=m\left(  s+l+1-i\right)  \right\}  $
for $2\leq i\leq l+1$. By Lemma \ref{maxppm} for each $i$ satisfying $1\leq
i\leq l-1,i\neq k$ there is $u_{i}$ such that $G_{i}=I_{u_{i}}$ and $u_{i}\neq
k,k+1$.

If $k=0$ then $G_{0}=\left[  1,\rho-1\right]  $ since all other cells are of
cardinality $\mu$. For each $G_{i}$ with $i\geq1$ the rank of the first
coordinate is $\rho+(i-1)\mu$, that is $r\left(  \beta,\rho+\left(
u_{i}-1\right)  \mu+1\right)  =\rho+\left(  i-1\right)  \mu$. Then by equation
\ref{keyq}
\[
\rho+\left(  i-1\right)  \mu-\left(  \rho+\left(  u_{i}-1\right)
\mu+1\right)  =\left(  \mu+1\right)  \left(  \left(  s+l+1-u_{i}\right)
-\left(  s+l+1-i\right)  \right)  ,
\]
thus $u_{i}=i+1$, for $1\leq i\leq l-1$. If $k>0$ we have shown $\#G_{0}=\rho$
and $\#G_{k}=\mu-1$. If $\rho\leq\mu-1$ then neither $G_{0}$ nor $G_{k}$ can
meet $I_{0}$ and another cell, by Lemma \ref{maxppm}. If additionally
$\rho<\mu-1$ then $G_{0}=I_{0}$ and $G_{k}=\left[  \rho+\left(  k-1\right)
\mu+1,\rho+k\mu-1\right]  $. If $\rho=\mu-1$ then it is not possible for
$G_{k}=I_{0}$ because then $r\left(  \beta,1\right)  =\rho+\left(  k-1\right)
\mu+1$ and equation \ref{keyq} yields $\rho+\left(  k-1\right)  \mu=\left(
\mu+1\right)  k$, that is $k+1=0$. As before, $G_{0}=I_{0}$ and $G_{k}%
=I_{k}\backslash\left\{  \rho+k\mu\right\}  $. If $\rho=\mu$ then $G_{k}%
=I_{k}\backslash\left\{  \rho+k\mu\right\}  $ and $G_{0}=I_{u_{0}}$ for some
$u_{0}$. The needed ranks for $\beta$ are $r\left(  \beta,\rho+\left(
u_{i}-1\right)  \mu+1\right)  =\rho+\left(  i-1\right)  \mu+1$ if $i<k$ and
$=\rho+\left(  i-1\right)  \mu$ if $k<i\leq l-1$. Similarly to the case $k=0$
this implies that $u_{i}=i$ for $i<k$ and $u_{i}=i+1$ for $k<i\leq l-1$.
\end{proof}

It remains to show that $\ell\left(  \beta\right)  \leq L$ implies
$\beta=\lambda$.

\begin{lemma}
Suppose $\beta_{L+1}=0$, then $\beta=\lambda$.
\end{lemma}

\begin{proof}
The hypothesis implies $\beta_{i}=0$ for all $i>L$ (the rank equation showed
$r\left(  \beta,L+1\right)  =L+1$ thus $i>L+1$ implies $\beta_{i}=0$). The
condition $\lambda\succeq\beta^{+}$ implies $\beta_{L}^{+}\geq\lambda
_{L}=m\left(  s+1\right)  $ (since $\left|  \lambda\right|  -\lambda_{L}%
\geq\left|  \beta\right|  -\beta_{L}^{+}$). For $1\leq i\leq l+1$ let
$G_{i}=\left\{  j:\beta_{j}=m\left(  s+l+1-i\right)  \right\}  $ and
$M_{i}=\#G_{i}$. Firstly let $m>1$, then $\beta_{\rho+k\mu}=m\left(
s+l+1-j_{0}\right)  -1$ for some $j_{0}$ in $0\leq j\leq l-1$, thus $\rho
+l\mu=\sum_{i=0}^{l}M_{i}+1$ and $\rho-M_{0}+\sum_{i=1}^{l}\left(  \mu
-M_{i}\right)  =1$. Also $M_{0}\leq\rho$ because $\lambda\succeq\beta^{+} $
and $M_{i}\leq\mu$ for $1\leq i\leq l$ by Lemma \ref{maxppm}. Hence either
$M_{0}=\rho-1$ and $M_{i}=\mu$ for $1\leq i\leq l$, or for some $j>0 $
$M_{j}=\mu-1$ and $M_{0}=\rho,M_{i}=\mu$ for $1\leq i\leq l,i\neq j. $ Now
\begin{align*}
\left|  \beta\right|   &  =\rho m\left(  s+l+1\right)  +\mu m\sum_{i=1}%
^{l}\left(  s+l+1-i\right)  -1\\
&  =M_{0}m\left(  s+l+1\right)  +m\sum_{i=1}^{l}M_{i}\left(  s+l+1-i\right)
+m\left(  s+l+1-j_{0}\right)  -1,
\end{align*}
and substituting the known values for $M_{i}$ we obtain $j_{0}=0$ if
$M_{0}=\rho-1$ else $j_{0}=j$. Then $r\left(  \beta,\rho+k\mu\right)
=\rho+j_{0}\mu$ and the rank equation at $\rho+k\mu$ yields $\left(
j_{0}-k\right)  \mu=\left(  \mu+1\right)  \left(  j_{0}-k\right)  $ and so
$j_{0}=k$, that is, $\beta_{\rho+k\mu}=\lambda_{\rho+k\mu}$. Similarly to the
previous lemma let $G_{i}=\left\{  j:\beta_{j}=m\left(  s+l+1-i\right)
\right\}  =I_{u_{i}}$ for $0\leq i\leq l,i\neq k$ and some $u_{i}\neq k$,
treating the special cases $\rho<\mu-1,\,\rho=\mu-1$ and $\rho=\mu$ as before.
Again $r\left(  \beta,\rho+\left(  i-1\right)  \mu+1\right)  =\rho+\left(
u_{i}-1\right)  \mu+1$ and the rank equation shows $u_{i}=i$. Also
$G_{k}=I_{k}\backslash\left\{  \rho+k\mu\right\}  $. Thus $\beta=\lambda.$

Secondly let $m=1$. Then $\lambda_{i}=s+l-k$ for $\rho+k\mu\leq i\leq
\rho+\left(  k+1\right)  \mu$. Also $\rho\geq M_{0}$ because $\lambda
\succeq\beta^{+}$. There are two equations involving $M_{i}$:
\begin{align}
\left(  \rho-M_{0}\right)  +\sum_{i=1}^{l}\left(  \mu-M_{i}\right)   &
=0,\label{msum1}\\
\left(  \rho-M_{0}\right)  \left(  s+l+1\right)  +\sum_{i=1}^{l}\left(
\mu-M_{i}\right)  \left(  s+l+1-i\right)   &  =1. \label{msum2}%
\end{align}
Equation \ref{msum2} shows that $\rho\geq M_{0}$ and $\mu\geq M_{i}$ for all
$i$ is impossible, hence there is at least one value, say $M_{j}$, such that
$M_{j}>\mu$. By the maximum principle $M_{j}=\mu+1$ and $M_{i}\leq\mu$ for all
$i\neq j$. Substituting these conditions in equation \ref{msum1} shows that
for some $j_{0}$, $M_{j_{0}}=\mu-1$ ($M_{0}=\rho-1$ if $j_{0}=0$) and
$M_{i}=\mu$ for all $i\neq j_{0},j$, and $M_{0}=\rho$ unless $j_{0}=0$.
Substitute these values in equation \ref{msum2} to obtain $j_{0}=j-1$. By
Lemma \ref{maxpp1} $G_{j}=\left[  \rho+k\mu,\rho+\left(  k+1\right)
\mu\right]  $, also $r\left(  \beta,\rho+k\mu\right)  =\rho+\left(
j-1\right)  \mu$. Then the rank equation shows $\left(  j-1-k\right)
\mu=\left(  \mu+1\right)  \left(  \left(  s+l-k\right)  -\left(
s+l+1-j\right)  \right)  =\left(  \mu+1\right)  \left(  j-1-k\right)  $ and
thus $j=k+1$. Similarly to the previous arguments, for each $i\neq k,k+1$
there exist $u_{i}$ such that $G_{i}=I_{u_{i}}$. Since $M_{k}+M_{k+1}=2\mu$
(or $\rho+\mu$ if $k=0$) we have $r\left(  \beta,\rho+\left(  u_{i}-1\right)
\mu+1\right)  =\rho+\left(  i-1\right)  \mu+1$ and the rank equation shows
$u_{i}=i$. This accounts for all of $\left[  1,L\right]  $ except for $\left[
1,\rho\right]  $ and $\left[  \rho+\left(  k-1\right)  \mu+1,\rho
+k\mu-1\right]  $. There are several cases for $\rho$: if $k=0$ then
$G_{0}=\left[  1,\rho-1\right]  $ by elimination; if $k\geq1$ and $\rho=\mu$
then $G_{0}=I_{u_{0}}$ and the rank equation shows $u_{0}=0$, and
$G_{k}=\left[  \rho+\left(  k-1\right)  \mu+1,\rho+k\mu-1\right]  $; if
$k\geq1$ and $\rho=\mu-1$ then by Lemma \ref{maxpp1} $G_{0}$ can not meet both
$I_{0}$ and $I_{k}$ thus either $G_{0}=I_{0}$ or $G_{0}=I_{k}\backslash
\left\{  \rho+k\mu\right\}  $ and the rank equation implies the latter can not
happen; if $k\geq1$ and $\rho<\mu-1$ then by the same Lemma $G_{k}%
=I_{k}\backslash\left\{  \rho+k\mu\right\}  $, forcing $G_{0}=I_{0}$. Thus
$\beta=\lambda.$
\end{proof}

The lemmas together provide the proofs of the following theorems.

\begin{theorem}
Let $\lambda=\Lambda\left(  \mu,s,l,\rho,m\right)  $ and $\lambda^{\left(
k\right)  }=\lambda-\varepsilon\left(  \rho+k\mu\right)  $ for $0\leq k\leq
l-1$. Then there exists a unique $\beta$ so that $\left(  \lambda^{\left(
k\right)  },\beta\right)  $ is $\left(  -\frac{m}{\mu+1}\right)  $-critical
and $\ell\left(  \beta\right)  =N+l-k>N$, where $N=\left(  s+l+1\right)
\mu+s+\rho$.
\end{theorem}

\begin{theorem}
Let $\lambda=\Lambda\left(  \mu,s,l,\rho,m\right)  $ and $N=\left(
s+l+1\right)  \mu+s+\rho$ then $\zeta_{\lambda}^{x}$ is a singular polynomial
for $S_{N}$ with singular value $-\frac{m}{\mu+1}$.
\end{theorem}

In the next section we study the irreducible representation associated to
$\zeta_{\lambda}^{x}$, in particular, an explicit basis for the span of its
$S_{N}$-orbit.

\section{Associated $S_{N}$-modules}

Using Murphy's construction \cite{Mu} of Young's seminormal representations we
can give a complete description of the $S_{N}$-orbit of $\zeta_{\lambda}^{x}.$
From the formula (valid for all $\kappa$ and for all polynomials $f$)
\[
\sum_{i=1}^{N}x_{i}\mathcal{D}_{i}f\left(  x\right)  =\sum_{i=1}^{N}x_{i}%
\frac{\partial}{\partial x_{i}}f\left(  x\right)  +\kappa\sum_{1\leq i<j\leq
N}\left(  f\left(  x\right)  -f\left(  x\left(  i,j\right)  \right)  \right)
\]
we note that a homogeneous singular polynomial $f$ must satisfy $\left(  \deg
f\right)  f=-\kappa\omega f$ where $\omega=\sum_{1\leq i<j\leq N}\left(
1-\left(  i,j\right)  \right)  $. But $\omega$ is in the center of
$\mathbb{Z}S_{N}$ and the eigenvalues for any isotype are known (Young's
formula). Indeed for any node $\left(  i,j\right)  $ in the Ferrers diagram of
a partition $\tau$ (with $\left|  \tau\right|  =N$), the content is defined to
be $c\left(  \left(  i,j\right)  \right)  =j-i$, then $\omega f=\left(
\binom{N}{2}-\sum_{\left(  i,j\right)  \in\tau}c\left(  \left(  i,j\right)
\right)  \right)  f$ whenever $f$ is of isotype $\tau$. Denote the eigenvalue
by $\tau\left(  \omega\right)  $, then $\tau\left(  \omega\right)  =\binom
{N}{2}-\frac{1}{2}\sum\limits_{i=1}^{\ell\left(  \tau\right)  }\tau_{i}\left(
\tau_{i}+1-2i\right)  $. As a function on partitions the eigenvalue is
strictly decreasing with respect to the dominance order.

\begin{lemma}
Suppose $\sigma,\tau\in\mathbb{N}_{0}^{N,P},\left|  \sigma\right|  =\left|
\tau\right|  $ and $\sigma\prec\tau$ then $\sum\limits_{\left(  i,j\right)
\in\sigma}c\left(  \left(  i,j\right)  \right)  <\sum\limits_{\left(
i,j\right)  \in\tau}c\left(  \left(  i,j\right)  \right)  $, and
$\sigma\left(  \omega\right)  >\tau\left(  \omega\right)  $.
\end{lemma}

\begin{proof}
By the theorems (1.15) and (1.16) in Macdonald \cite[p.9]{Ma} it suffices to
prove the inequality for $\tau=\sigma+\varepsilon\left(  i\right)
-\varepsilon\left(  j\right)  $ with $i<j$ (this is a ``raising operator'').
Then $\sum\limits_{\left(  i,j\right)  \in\tau}c\left(  \left(  i,j\right)
\right)  =\sum\limits_{\left(  i,j\right)  \in\sigma}c\left(  \left(
i,j\right)  \right)  +\left(  \sigma_{i}-\sigma_{j}\right)  +\left(
j+1-i\right)  $.
\end{proof}

Recall the singular polynomials $\zeta_{\lambda}^{x}$ associated to two-part
partitions $\tau=\left(  \mu,N-\mu\right)  $ with $\lambda=\left(  m^{N-\mu
}\right)  $ and $\gcd\left(  m,\mu+1\right)  <\frac{\mu+1}{N-\mu}$; then
$\tau\left(  \omega\right)  =\left(  \mu+1\right)  \left(  N-\mu\right)  $ and
$\deg\zeta_{\lambda}^{x}=m\left(  N-\mu\right)  $. For $\tau=\left(  s\left(
\mu+1\right)  +\mu,\mu^{l},\rho\right)  $ we find
\begin{align*}
\tau\left(  \omega\right)   &  =\left(  \mu+1\right)  \left(  \rho\left(
s+l+1\right)  +\frac{1}{2}\mu l\left(  l+2s+1\right)  \right) \\
&  =\frac{\mu+1}{m}\left|  \Lambda\left(  \mu,s,l,\rho,m\right)  \right|  .
\end{align*}

\begin{theorem}
For $\lambda=\Lambda\left(  \mu,s,l,\rho,m\right)  $ and $\kappa=-\frac{m}%
{\mu+1}$ the singular polynomial $\zeta_{\lambda}^{x}$ on $\mathbb{R}^{N} $ is
of isotype $\tau=\left(  s\left(  \mu+1\right)  +\mu,\mu^{l},\rho\right)  $
($\left|  \tau\right|  =N$).
\end{theorem}

\begin{proof}
For any $\zeta_{\sigma}^{x}$ with $\sigma\in\mathbb{N}_{0}^{N}$ if $\sigma
_{i}=\sigma_{i+1}$ for some $i$ then $\left(  i,i+1\right)  \zeta_{\sigma}%
^{x}=\zeta_{\sigma}^{x}$. Thus $\zeta_{\lambda}^{x}$ is invariant under
$S_{\left[  1,\rho\right]  }\times\prod_{j=1}^{l}S_{\left[  \rho+\left(
j-1\right)  \mu+1,\rho+j\mu\right]  }\times S_{\left[  \rho+l\mu+1,N\right]
}$, and this group is conjugate to $S_{\tau}$ (the direct product $\prod
_{i}S_{\tau_{i}}$). Thus $E=\mathrm{span}_{\mathbb{Q}}\left\{  w\zeta
_{\lambda}^{x}:w\in S_{N}\right\}  $ is isomorphic to a submodule of the
representation of $S_{N}$ induced up from $1_{S_{\tau}},$ the identity
representation of $S_{\tau}$. By a classical theorem (see Macdonald
\cite[p.115]{Ma}) this decomposes as a direct sum with one component of
isotype $\tau$ and all other components of isotypes $\sigma$ with $\sigma
\succ\tau$. Any $f\in E$ is singular and $f$ can not be of isotype
$\sigma\succ\tau$ because $\deg f=\allowbreak\left|  \lambda\right|
=\allowbreak\frac{m}{\mu+1}\tau\left(  \omega\right)  >\frac{m}{\mu+1}%
\sigma\left(  \omega\right)  $ by the Lemma.
\end{proof}

The same method proves the following.

\begin{proposition}
For $\frac{N}{2}\leq\mu<N$, $\gcd\left(  m,\mu+1\right)  <\frac{\mu+1}{N-\mu}%
$, and $\lambda=\left(  m^{N-\mu}\right)  $ the singular polynomial
$\zeta_{\lambda}^{x}$ on $\mathbb{R}^{N}$ for $\kappa=-\frac{m}{\mu+1}$ is of
isotype $\left(  \mu,N-\mu\right)  $.
\end{proposition}

We turn to the application of Murphy's results. For any given isotype he
determined the eigenvalues and eigenvectors of the commuting operators
\newline$\left\{  \sum\limits_{j=1}^{i-1}\left(  i,j\right)  :2\leq i\leq
N\right\}  $ (Jucys-Murphy elements). However the results have to be read in
reverse in a certain sense.

\begin{proposition}
Suppose $f$ is a singular polynomial for $\kappa=\kappa_{0}\in\mathbb{Q}$ and
$1\leq i\leq N$, then $\mathcal{U}_{i}f=f+\kappa_{0}\sum_{j=i+1}^{N}\left(
i,j\right)  f$.
\end{proposition}

\begin{proof}
We have the commutation $\mathcal{D}_{i}\left(  x_{i}f\right)  =x_{i}%
\mathcal{D}_{i}f+f+\kappa\sum_{j\neq i}\left(  i.j\right)  f$. Now set
$\kappa=\kappa_{0}$ and note that $\mathcal{U}_{i}f=\mathcal{D}_{i}\left(
x_{i}f\right)  -\kappa\sum_{j<i}\left(  i,j\right)  f$.
\end{proof}

Denote the Murphy elements $\omega_{i}=\sum\limits_{j=N-i+2}^{N}\left(
N+1-i,j\right)  $ for $2\leq i\leq N$ and let $\omega_{1}=0$ (as a
transformation); then $\mathcal{U}_{i}f=f+\kappa_{0}\omega_{N+1-i}f$ for
singular polynomials. Suppose that $\zeta_{\alpha}^{x}$ is singular for
$\kappa=\kappa_{0}$ and $\alpha=w\lambda$, some $w\in S_{N}$ (recalling that
$\mathcal{U}_{i}\zeta_{\alpha}^{x}=\left(  \left(  N-r\left(  \alpha,i\right)
\right)  \kappa+\alpha_{i}+1\right)  \zeta_{\alpha}^{x}$), then $\omega
_{N+1-i}\zeta_{\alpha}^{x}=\left(  N-r\left(  \alpha,i\right)  +\frac
{\alpha_{i}}{\kappa_{0}}\right)  \zeta_{\alpha}^{x}$. A \textit{standard Young
tableau} (SYT) of shape $\tau$ is a one-to-one assignment of the numbers
$\left\{  1,\ldots,N\right\}  $ to the nodes of the Ferrers diagram so that
the entries increase in each row and in each column. Let $\eta_{i}\left(
T\right)  $ be the content of the node containing the value $i$, $1\leq i\leq
N$. Murphy constructed a basis $\left\{  f_{T}:T\text{ is an SYT of shape
}\tau\right\}  $ for the irreducible representation of isotype $\tau$ and
$\omega_{i}f_{T}=\eta_{i}\left(  T\right)  f_{T}$ for each $i$ and $T$. There
is an order on SYT's of given shape (for details see \cite[p.288]{Mu}) and the
maximum SYT in this order, denoted by $T_{0}$, is produced by entering the
numbers $1,2,\ldots,N$ row by row (the first row is $1,\ldots,\tau_{1}$, the
second is $\tau_{1}+1,\ldots,\tau_{1}+\tau_{2}$ and so forth).

\begin{definition}
Suppose $T$ is an SYT\ of shape $\tau$, with $\tau\in\mathbb{N}_{0}^{N,P}$ and
$\left|  \tau\right|  =N$, then let $rw\left(  i,T\right)  ,cm\left(
i,T\right)  $ denote the row and column respectively of the node of $T$
containing $i$, for $1\leq i\leq N$. Let $t\left(  i,\tau\right)  $ (or
$t\left(  i\right)  )=\left(  rw\left(  i,T_{0}\right)  ,cm\left(
i,T_{0}\right)  \right)  $, considered as a labeling of the nodes in the
diagram of $\tau$.
\end{definition}

In this notation $\eta_{i}\left(  T\right)  =cm\left(  i,T\right)  -rw\left(
i,T\right)  $.

\begin{proposition}
\label{lb2cont}Let $\lambda=\Lambda\left(  \mu,s,l,\rho,m\right)  $
(hypotheses as in Definition \ref{lambda}) then $N-k+\frac{\lambda_{k}}%
{\kappa_{0}}=c\left(  t\left(  N+1-k\right)  \right)  =\eta_{N+1-k}\left(
T_{0}\right)  $ for $1\leq k\leq N$
\end{proposition}

\begin{proof}
Let $c_{k}=N-k+\frac{\lambda_{k}}{\kappa_{0}}=N-k-\left(  \mu+1\right)
\frac{\lambda_{k}}{m}$ for $1\leq k\leq N$, then for $k=\rho+1-j$ and $1\leq
j\leq\rho$ we have $c_{k}=j-\left(  l+2\right)  $, for $k=\rho+\left(
l+1-i\right)  \mu+1-j$ with $1\leq i\leq l$ and $1\leq j\leq\mu$ we have
$c_{k}=j-\left(  i+1\right)  $, and finally for $k=N+1-j$ with $1\leq j\leq
N-\ell\left(  \lambda\right)  =s\left(  \mu+1\right)  +\mu=\tau_{1}$ we have
$c_{k}=j-1$. Thus $c_{k}=c\left(  t\left(  N+1-k\right)  \right)  .$
\end{proof}

For a partition $\lambda\in\mathbb{N}_{0}^{N,P}$ say that $w\in S_{N}$ is
$\lambda$-rank-preserving if $\lambda_{i}=\lambda_{i+1}$ implies $w\left(
i\right)  <w\left(  i+1\right)  $ for $1\leq i<N$. In general $\left(
w\lambda\right)  _{w\left(  i\right)  }=\lambda_{i}$, so this property implies
$r\left(  w\lambda,i\right)  =r\left(  \lambda,w^{-1}\left(  i\right)
\right)  =w^{-1}\left(  i\right)  $ for $1\leq i\leq N$ and
\[
\mathcal{U}_{i}\zeta_{w\lambda}^{x}=\left(  \left(  N-w^{-1}\left(  i\right)
\right)  \kappa+\lambda_{w^{-1}\left(  i\right)  }+1\right)  \zeta_{w\lambda
}^{x}%
\]
for generic $\kappa$. In particular, if $\lambda=\Lambda\left(  \mu
,s,l,\rho,m\right)  $, $w$ is $\lambda$-rank-preserving and $\zeta_{w\lambda
}^{x}$ is singular (for $\kappa=\kappa_{0}$) then
\[
\omega_{N+1-i}\zeta_{w\lambda}^{x}=c\left(  t\left(  N+1-w^{-1}\left(
i\right)  \right)  \right)  \zeta_{w\lambda}^{x}.
\]
Let $w_{0}$ be the \textquotedblleft reversing\textquotedblright\ (longest)
element of $S_{N}$, that is $w_{0}\left(  i\right)  =N+1-i$ for $1\leq i\leq
N$ (note $w_{0}^{-1}=w_{0}$). Thus $\omega_{i}\zeta_{w\lambda}^{x}=c\left(
t\left(  w_{0}ww_{0}\left(  i\right)  \right)  \right)  \zeta_{w\lambda}^{x}$.
On the other hand suppose $u\in S_{N}$ and the action of $u$ on $T_{0}$
produces an SYT denoted by $T$ ($u$ acts on the entries of $T_{0}$), then the
node $t\left(  i\right)  $ contains $u\left(  i\right)  $. Thus $\omega
_{i}f_{T}=\eta_{i}\left(  T\right)  f_{T}=c\left(  t\left(  u^{-1}\left(
i\right)  \right)  \right)  f_{T}$ and $f_{T}$ has the same respective
eigenvalues for $\left\{  \omega_{i}:1\leq i\leq N\right\}  $ as
$\zeta_{w\lambda}^{x}$ for $w=w_{0}^{-1}uw_{0}$ provided that $w$ is $\lambda
$-rank-preserving. But this is a consequence of $T$ being an SYT ($\lambda
_{i}=\lambda_{i+1}$ implies $N-i$ and $N+1-i$ are in the same row of $T_{0}$
thus $u\left(  N-i\right)  $ and $u\left(  N-i+1\right)  $ are in the same row
of $T,$with $u\left(  N-i\right)  <u\left(  N-i+1\right)  $, that is,
$w\left(  i\right)  <w\left(  i+1\right)  $). Further it is easy to describe
$w\lambda$ corresponding to a given SYT $T$: let $\gamma_{1}=0$ and
$\gamma_{j}=m\left(  s+j-1\right)  $ for $2\leq j\leq l+2$, then for any $i$
(with $1\leq i\leq N$) let $j=rw\left(  i,T\right)  $ and $\left(
w\lambda\right)  _{N+1-i}=\gamma_{j}$. This shows that the possible $w\lambda$
corresponding to SYT's are exactly the reverse lattice permutations of
$\lambda$. A \textit{reverse lattice permutation} $w\lambda$ of $\lambda$ is
defined by the property that every right substring $\left(  w\lambda\right)
_{N+1-j}\left(  w\lambda\right)  _{N+2-j}\ldots\left(  w\lambda\right)  _{N}$
(for $1\leq j\leq N$) has at least as many entries of $\gamma_{i}$ as of
$\gamma_{i+1}$ for each $i$. The set of corresponding $w\in S_{N}$ serves as
an index set, namely $E_{\tau}=\left\{  w\in S_{N}:w_{0}ww_{0}T_{0}\text{ is
an SYT}\right\}  $.

We show that $\left\{  \zeta_{w\lambda}^{x}:w\in E_{\tau}\right\}  $ is a
basis for the $S_{N}$-module (isotype $\tau$) generated by $\zeta_{\lambda
}^{x}$.

\begin{theorem}
\label{Murf}Let $w\in E_{\tau}$ then $\zeta_{w\lambda}^{x}$ is singular (for
$\kappa=\kappa_{0}$ on $\mathbb{R}^{N}$) and $\left\{  \zeta_{w\lambda}%
^{x}:w\in E_{\tau}\right\}  $ is a basis for $span_{\mathbb{Q}}\left\{
w\zeta_{\lambda}^{x}:w\in S_{N}\right\}  $, on which $S_{N}$ acts by Young's
seminormal representation, where $\zeta_{w\lambda}^{x}$ corresponds to $f_{T}$
with $T=w_{0}ww_{0}T_{0}$.
\end{theorem}

\begin{proof}
By Proposition \ref{z2sz} if $\zeta_{w\lambda}^{x}$ has no pole at $\kappa
_{0}$, for some $w\in S_{N}$, and $a=\kappa\left(  \kappa\left(  r\left(
w\lambda,i+1\right)  -r\left(  w\lambda,i\right)  \right)  +\lambda
_{w^{-1}\left(  i\right)  }-\lambda_{w^{-1}\left(  i+1\right)  }\right)
^{-1}$ does not evaluate to $\pm1$ at $\kappa=\kappa_{0},$for some $i$ with
$\lambda_{w^{-1}\left(  i\right)  }>\lambda_{w^{-1}\left(  i+1\right)  }$ then
$\zeta_{\left(  i,i+1\right)  w\lambda}^{x}$ does not have a pole at
$\kappa_{0}$ (the formula is $\left(  i,i+1\right)  \zeta_{w\lambda}%
^{x}=a\zeta_{w\lambda}^{x}+\left(  1-a^{2}\right)  \zeta_{\left(
i,i+1\right)  w\lambda} $). By Proposition \ref{lb2cont} $a=\left(  c\left(
t\left(  N+1-w^{-1}\left(  i\right)  \right)  \right)  -c\left(  t\left(
N+1-w^{-1}\left(  i+1\right)  \right)  \right)  \right)  ^{-1}$ at
$\kappa=\kappa_{0}$. Each SYT $T$ of shape $\tau$ is the result of a (finite)
sequence $\ \left\{  \left(  i_{j},i_{j}+1\right)  :1\leq j\leq n\right\}  $
of adjacent transpositions applied to $T_{0}$, such that if $T_{j}=\left(
i_{j},i_{j}+1\right)  T_{j-1}$ then $rw\left(  i_{j},T_{j-1}\right)
<rw\left(  i_{j}+1,T_{j-1}\right)  $. This also implies that $T_{j}$ is lower
in the order on tableaux as used in \cite{Mu}.

For any SYT $T$ there are four possibilities for the locations of $i,i+1$ and
Murphy \cite[p.292]{Mu} derived the expansion of $\left(  i,i+1\right)  f_{T}$
in each case: if $rw\left(  i,T\right)  =rw\left(  i+1,T\right)  $ then
$\left(  i,i+1\right)  f_{T}=f_{T}$, if $cm\left(  i,T\right)  =cm\left(
i+1,T\right)  $ then $\left(  i,i+1\right)  f_{T}=-f_{T}$ and if $rw\left(
i,T\right)  <rw\left(  i+1,T\right)  $ then $\left(  i,i+1\right)
f_{T}=af_{T}+\left(  1-a^{2}\right)  f_{\left(  i,i+1\right)  T}$ where
$a=\left(  \eta_{i}\left(  T\right)  -\eta_{i+1}\left(  T\right)  \right)
^{-1}$ (the fourth case, $rw\left(  i,T\right)  >rw\left(  i+1,T\right)  $
follows from the previous by interchanging $T$ and $\left(  i,i+1\right)  T$;
also $rw\left(  i,T\right)  <rw\left(  i+1,T\right)  $ implies $cm\left(
i,T\right)  >cm\left(  i+1,T\right)  $ thus $0<a\leq\frac{1}{2}$). As remarked
before, if $w\lambda$ corresponds to an SYT $T$ with $rw\left(  i,T\right)
<rw\left(  i+1,T\right)  $ then $\left(  w\lambda\right)  _{N+1-i}<\left(
w\lambda\right)  _{N-i}$. Let $\beta=\left(  N-i,N-i+1\right)  w\lambda$, then
$\zeta_{\beta}^{x}=\left(  1-a^{2}\right)  ^{-1}\left(  \left(
N-i,N-i+1\right)  \zeta_{w\lambda}^{x}-a\zeta_{w\lambda}^{x}\right)  $ with
the same $a$ that appears in the expression for $f_{\left(  i,i+1\right)  T}$
in terms of $f_{T}$ (note $w_{0}\left(  i,i+1\right)  w_{0}=\left(
N-i,N-i+1\right)  $). Since $f_{T_{0}}$ has the same eigenvalues for $\left\{
\omega_{i}\right\}  $ as $\zeta_{\lambda}^{x}$ this argument used inductively
(on the number of adjacent transpositions linking $T_{0}$ to $T$) shows that
$\left\{  \zeta_{w\lambda}^{x}:w\in E_{\tau}\right\}  $ transforms according
to the seminormal representation, for the isotype $\tau$. Again suppose
$\zeta_{w\lambda}^{x}$ corresponds to the SYT $T$ (that is $T=w_{0}ww_{0}%
T_{0}$); if $rw\left(  i,T\right)  =rw\left(  i+1,T\right)  $ then $\left(
w\lambda\right)  _{N+1-i}=\left(  w\lambda\right)  _{N-i}$ and $\zeta
_{w\lambda}^{x}$ is invariant under $\left(  N-i,N-i+1\right)  ,$while if
$cm\left(  i,T\right)  =cm\left(  i+1,T\right)  $ then $\eta_{i}\left(
T\right)  -\eta_{i+1}\left(  T\right)  =1$ and the equation $\left(
N-i,N-i+1\right)  \zeta_{w\lambda}^{x}=-\zeta_{w\lambda}^{x}$ is a consequence
of the fact that $S_{N}$ acts on the basis $\left\{  \zeta_{w\lambda}^{x}:w\in
E_{\tau}\right\}  $ just as on $\left\{  f_{T}\right\}  $.
\end{proof}

The concept of reverse lattice permutations of $\lambda$ provides a concise
labeling of the singular polynomials of isotype $\tau$.

\section{Conclusion}

Here is a description of how to find the isotype $\tau$ and label $\lambda$
for the singular value $\kappa=-\frac{m}{n}$, given a pair $\left(
m,n\right)  $ with $2\leq n\leq N,m\geq1$ and $\frac{m}{n}\notin\mathbb{N}$.
Let $d=\gcd\left(  m,n\right)  ,\,m_{1}=\frac{m}{d},\,n_{1}=\frac{n}{d}$ (by
hypothesis $n_{1}\geq2$), then let $l=\left\lceil \dfrac{N+1-n}{n_{1}%
-1}\right\rceil -1$ (the ceiling function), $\rho=\left(  N+1-n\right)
-l\left(  n_{1}-1\right)  $. If $l=0$ then $\tau=\left(  n-1,N+1-n\right)  $
and $\lambda=\left(  m^{N+1-n},0^{n-1}\right)  .$ If $l\geq1$ then
$\tau=\left(  n-1,\left(  n_{1}-1\right)  ^{l},\rho\right)  $, $\mu
=n_{1}-1,\,s=d-1$ and $\lambda=\Lambda\left(  n_{1}-1,d-1,l,\rho,m_{1}\right)
$, that is, $\lambda=\left(  \left(  m+lm_{1}\right)  ^{\rho},\left(
m+\left(  l-1\right)  m_{1}\right)  ^{\mu},\ldots,m^{\mu},0^{n-1}\right)  $.
Note that the first part of $\tau$ is always $n-1$ (and $\lambda$ ends in
$n-1$ zeros).

The rational Cherednik algebra $\mathbf{A}$ was investigated by Berest,
Chmutova, Etingof, Ginzburg, Guay, Opdam and Rouquier in a series of papers
\cite{BEG1},\cite{BEG2},\cite{CE},\cite{EG},\cite{GGOR}. Here we consider the
faithful representation of $\mathbf{A}$ as the algebra generated by
\newline$\left\{  \mathcal{D}_{i},x_{i}:1\leq i\leq N\right\}  \cup S_{N}$ of
operators on polynomials on $\mathbb{R}^{N}$ (where $x_{i}$ denotes the
multiplication operator). Suppose $\kappa=\kappa_{0}$ is a singular value and
$\tau,\lambda$ are defined as in the previous section and used in Theorem
\ref{Murf}, and let $M_{\tau}=\mathrm{span}_{\mathbb{P}}\left\{
\zeta_{w\lambda}^{x}:w\in E_{\tau}\right\}  $ where $\mathbb{P}$ denotes
$\mathbb{Q}\left[  x_{1},\ldots,x_{N}\right]  ,$ the polynomials on
$\mathbb{R}^{N}$. Then $M_{\tau}$ is a module for the Cherednik algebra
$\mathbf{A}$ specialized to $\kappa=\kappa_{0}$; clearly $M_{\tau}$ is closed
under multiplication by polynomials and the action of $S_{N}.$ It is closed
under $\left\{  \mathcal{D}_{i}:1\leq i\leq N\right\}  $, indeed suppose $p$
is a polynomial and $g\in\mathrm{span}_{\mathbb{Q}}\left\{  \zeta_{w\lambda
}^{x}:w\in E_{\tau}\right\}  $ then by the product rule $\mathcal{D}%
_{i}\left(  pg\right)  =p\mathcal{D}_{i}g+g\dfrac{\partial}{\partial x_{i}%
}p+\kappa\sum\limits_{j\neq i}\left(  \left(  i,j\right)  g\right)
\dfrac{p\left(  x\right)  -\left(  i,j\right)  p\left(  x\right)  }%
{x_{i}-x_{j}}\in M_{\tau}$ when $\kappa=\kappa_{0}$ and $\mathcal{D}_{i}g=0$.

It is a plausible conjecture that we have found all the singular polynomials
for $S_{N}$ (perhaps to be settled in a later paper). The structure of
$\kappa_{0}$-critical pairs (see Definition \ref{critpair}) may be worth
further investigation, with a view to finding a general algorithm for their
construction, and maybe a uniqueness result in the case that $h\left(
\lambda,\kappa+1\right)  $ has a zero of multiplicity 1 at $\kappa=\kappa_{0}%
$. Such a result would simplify the argument used here.

\end{document}